\documentclass{gtpart}
\usepackage[english]{babel}
\usepackage{euscript}
\usepackage{amsfonts}
\usepackage{amsmath, amscd, amssymb}
\usepackage{xy}
\input xy
\xyoption{all} 
\usepackage{latexsym}

\def \M{\mathbf{MO_{R}}(W)}
\def \I{\EuScript{I}}
\def \K{\mathcal{K}}
\def \Symm{\K}

\def \R{\mathbb{R}}
\def \Z{\mathbb{Z}}

\def \J{\EuScript{J}}
\def \EB{\EuScript{B}}
\def \EC{\EuScript{C}}
\def \B{\mathbf{B}}
\def \RO{{R}}

\def \BDiff{\mathop{\mathrm{BDiff}}}

\def \Diff{\mathop{\mathrm{Diff}}}
\def \id{\mathrm{id}}
\def \jet{\mathrm{jet}}
\def \ie {i.e.}

\def \Exp{\mathrm{Exp}}

\def \Top{\mathbf{Top}}

\title{Bordism groups of solutions to differential relations}

\author{Rustam Sadykov}

\address{Department of Mathematics \\
         University of Toronto}

\email{sadykov@math.toronto.edu}


\volumenumber{}
\issuenumber{}
\publicationyear{}
\papernumber{}
\startpage{}
\endpage{}
\doi{}
\MR{}
\Zbl{}
\received{}
\revised{}
\accepted{}
\published{}
\publishedonline{}
\proposed{}
\seconded{}
\corresponding{}
\editor{}
\version{}

\newcommand{\colim}{\mathop{\mathrm{colim}}}

\newtheorem{theorem}{Theorem}[section]
\newtheorem*{thM}{Main Theorem}

\newtheorem{lemma}[theorem]{Lemma}
\newtheorem{corollary}[theorem]{Corollary}

\theoremstyle{definition}
\newtheorem{definition}[theorem]{Definition}
\newtheorem*{defB}{Definition of $\B_{\J}$}
\newtheorem{remark}[theorem]{Remark}

\newtheorem{example}[theorem]{Example}

\begin{document}

\begin{abstract}
In terms of category theory, the {\it Gromov homotopy principle} for a set valued functor $F$ asserts that the functor $F$ can be induced from a homotopy functor. Similarly, we say that the {\it bordism principle} for an abelian group valued functor $F$ holds if the functor $F$ can be induced from a (co)homology functor. 

We examine the bordism principle in the case of functors given by (co)bordism groups of maps with prescribed singularities. Our main result implies that if a family $\J$ of prescribed singularity types satisfies certain mild conditions, then there exists an infinite loop space $\Omega^{\infty}\B_{\J}$ such that for each smooth manifold $W$ the cobordism group of maps into $W$ with only $\J$-singularities is isomorphic to the group of homotopy classes of maps $[W, \Omega^{\infty}\B_{\J}]$. The spaces $\Omega^{\infty}\B_{\J}$ are relatively simple, which makes explicit computations possible even in the case where the dimension of the source manifold is bigger than the dimension of the target manifold.  

\end{abstract}
%
%

\maketitle

\section{Introduction}\label{s:1}

A smooth map $f\co V\to W$ of manifolds is said to be {\it singular} at a point $x\in V$ if the rank of the differential of $f$ at $x$ is less than the minimum of dimensions of $V$ and $W$. If $x$ is a {\it regular value} of $f$, i.e., if $f^{-1}(x)$ consists only of non-singular points, then $f^{-1}(x)$ is a smooth manifold of a {\it formal dimension} $d=\dim{V}-\dim{W}$. The integer $d$ is also called the {\it dimension} of $f$. In the present paper we are primarily interested in the case $d\ge 0$.

For a set of prescribed singularity types $\J$ of maps of a fixed dimension $d$, a smooth map $f$ of manifolds is said to be a {\it $\J$-map} if each singular point of $f$ is of type in $\J$. Similarly, we say that a cobordism of two $\J$-maps is a {\it $\J$-cobordism} if, as a map, it has only $\J$-singularities. The set of $\J$-cobordism classes of $\J$-maps of closed manifolds into a closed manifold $W$ leads to an abelian group $B(W; \J)$; in terms of representatives the group operation is given by taking the disjoint union of maps (see sections~\ref{s:5}, \ref{s:6}). 

Let us mention that in later sections we will define the group $B(W; \J)$ for an arbitrary (not necessarily closed) manifold $W$. At the moment, however, we consider only maps of closed manifolds, just to simplify the exposition.  

Carefully choosing the set $\J$ of singularities, one may derive cobordism groups related to various objects in geometry and topology. 

\begin{example}
Since a proper submersion is a smooth fiber bundle, the cobordism group of submersions is closely related to diffeomorphism groups of smooth manifolds.
It is also known to be related to various infinite loop spaces, moduli spaces of Riemann surfaces, the cobordism category \cite{GMTW}; and, in particular, to the Kahn-Priddy theorem and the standard Mumford conjecture, \cite{Fu}, \cite{MW}, \cite{GMTW} (see Remark~\ref{r:2}).
Similarly, $\J$-cobordism groups are related to singular cobordism categories~\cite{Sad3} and a certain version of $\J$-cobordisms is related  to the Deligne-Mumford compactification~\cite{EG}. The cobordism groups of {\it special generic maps} are related to dif\-feomor\-phism groups of spheres and exotic smooth structures on spheres~\cite{Sae1}, \cite{Sad1}.  The cobordism groups of {\it fold maps} are known to be related to stable homotopy groups of spheres and other interesting objects~\cite{Ch}, \cite{An1}, \cite{Kal2}. 
\end{example}

A priori $\J$-cobordism groups do not form generalized cohomology theories since, for example, $\J$-cobordism groups are not defined for topological spaces. In the current paper we propose a counterpart of $B(W; \J)$ that for a wide range of sets $\J$ can be used to compute $B(W; \J)$ in the same way as singular cohomology groups $H^n(W; \R)$ can be used to compute De Rham cohomology groups $H^n_{DR}(W)$.

%

\begin{definition}\label{def:1.1} 
Let $F$ be a contravariant functor from a category $\mathcal{C}$ to the category $\mathbf{Ab}$ of abelian groups. We say that $F$ satisfies the {\it  bordism principle} (or {\it  b-principle}) if there is a cohomology theory with functors $h^i\co  \mathbf{Top}\to \mathbf{Ab}$, indexed by $i\in \Z$, on the category $\mathbf{Top}$ of topological spaces and a covariant factor $\tau\co  \mathcal{C} \to \mathbf{Top}$ such that $F$ is naturally equivalent to $h^n\circ \tau$ for some $n$, i.e., there is a non-commutative diagram of functors
\[
\xymatrix{
  &  \mathcal{C} \ar[rr]^{\tau} \ar[dr]_{F} &     &  \mathbf{Top} \ar[dl]^{h^n} & \\
  &                                 & \mathbf{Ab}   &                     &
}
\] 
\end{definition}

On the category of smooth manifolds, for example, $H^n_{DR}$ satisfies the b-principle since it is naturally equivalent to $H^n\circ \tau$; here $\tau$ takes a smooth manifold onto the underlying topological space. 

In general, if holds, the b-principle allows us to replace an a priori incomprehensible functor by one that can be studied by means of the machinery of the cohomology theory.

\begin{remark} The b-principle is a bordism version of the {\it homotopy principle}, or {\it h-principle}. The classical Gromov definition of the h-principle is given in terms of so-called {\it jet bundles} \cite{Gr}. There are, however, important h-principle type theorems, e.g., Thurston h-principle for foliations \cite{Th1}, \cite{Th2} (see also \cite{EM0}, \cite{EM1}) that do not fit the
classical jet bundle setting. In general, the h-principle can be formulated in terms of category theory as above by respectively replacing $\mathbf{Top}$ and $\mathbf{Ab}$ by the category $\mathbf{Top}\downarrow B$  of topological spaces over a fixed space $B$ (see section~\ref{s:2}) and a category of sets, and requiring that $h^n$ be a homotopy functor rather than a term of a cohomology theory. 
\end{remark}

In the present paper we show that for a wide range of sets $\J$ the b-principle holds true for a functor counterpart of $B(*; \J)$. More importantly, we construct a cohomology theory extending $B(*; \J)$ whose spectrum $\B_{\J}$ is simple enough to make  explicit computations possible. 

\subsection*{Structure of the paper} The results of the paper are stated in section~\ref{s:1.1}. In section~\ref{s:2} we formulate the bordism and weak bordism principles. In sections~\ref{s:rel}-\ref{s5} and \ref{s:JVB}-\ref{h-principle} we recall and develop the language necessary to discuss (co)bordism groups in terms of jet spaces and differential equations and, more generally, differential relations. We will see that each set $\J$ of singularity types is associated with differential relations so that $\J$-maps can be defined to be solutions of associated differential relations. In sections \ref{s:5}-\ref{s:6} we define bordism and cobordism groups of solutions of differential relations. In sections
\ref{s:iso}-\ref{s:13} we prove the Main Theorem and its covariant counterpart. Final sections~\ref{s:app} and \ref{s:15} are devoted to applications. 


\subsection*{Acknowledgment} I gratefully acknowledge that Theorem~\ref{th:obst} was poin\-ted out to me by Andr\'{a}s Sz\H{u}cs. I am thankful to Osamu Saeki, Alex Dranishnikov and Yuli Rudyak for their time, attention and numerous suggestions while the paper was written, and to Yoshifumi Ando for discussions of the results of the paper. I am also grateful to the referee for many comments which lead to improvement in the exposition of results. 

The author has been supported by a Postdoctoral Fellowship and a Grant-in-Aid for Scientific Research of the Japan Society for the Promotion of Science, and a scholarship of Max Planck Institute of Mathematics.

\section{Results}\label{s:1.1}

\begin{defB} Let $p\co  EO_t\to BO_t$ be the universal vector $t$-bundle. Let $S_t=S_t(\J)$ denote the space of Taylor series $T(f)$ at $0\in \R^{t+d}$ of $\J$-maps $f\co  \R^{t+d} \to EO_t|b$,
with $b\in BO_t$, into the fibers $EO_t|b=p^{-1}(b)$ such that $f(0)$ is $0\in EO_t|b$. Then the map $\pi\co  S_t\to BO_t$ that takes $T(f)$ onto $b$ has a structure of a fiber bundle. The desired spectrum $\B_{\J}$ is defined to be the Thom spectrum with $(t+d)$-th term given by the Thom space $T\pi^*EO_t$ of the bundle $\pi^*EO_t$
over $S_t$.  
\end{defB}

Unless otherwise stated, we will always (tacitly) assume that the set $\J$ of (non-regular) singularity types is {\it open, $\K$-invariant}. The former means that a map close to a $\J$-map is also a $\J$-map, while the latter means that the set $\J$ has {\it sufficiently many symmetries}  (see section~\ref{s:rel}). In addition, unless otherwise stated, we will always (tacitly) assume that if $d\ge 0$, then $\J$ contains folds and the target $W$ of considered maps is of dimension $>1$.

\begin{thM} The cobordism group $B(W; \J)$ of $\J$-maps of dimension $d$ to a closed manifold $W$ is isomorphic to the set $[W, \Omega^{\infty}\mathbf{B}_{\J}]$ of homotopy classes of maps into the infinite loop space $\Omega^{\infty}\mathbf{B}_{\J}$ of $\mathbf{B}_{\J}$. 
\end{thM}

\begin{remark} Similar definitions and a theorem are valid for bordism groups of $\J$-maps and their functor counterparts (see section~\ref{s:5} and Theorem~\ref{t:iso}).
\end{remark}

\begin{remark}\label{r:2} In the case of maps of a general dimension $d$ with $\J=\emptyset$, i.e., in the case where $\J$-maps have only non-singular points, 
the statement of the Main Theorem is not true. Let us mention, however, that its version, which we do not consider in the present paper, holds true for $d=0$ (Kahn-Priddy theorem, \cite{Ch}, \cite{Fu}) and $d=2$ (Mumford conjecture, \cite{MW}).  
\end{remark}

\begin{remark} The omitted case of maps of dimension $d\ge 0$ into a manifold $W$ of dimension $1$ has been considered in~\cite{IS}, \cite{Ik}, \cite{Kal0}, \cite{Sae1}. For this case our approach does not apply because our argument implicitly (but essentially) uses the Eliashberg h-principle~\cite{El1}, \cite{El2} for fold maps which requires the condition $\dim W>1$.  
\end{remark}

The Main Theorem is obviously related to an old question in singularity theory on constructing a classifying space $B_\J$ such that for each closed manifold $W$ there is an isomorphism $B(W; \J)\approx [W, B_{\J}]$. The first general result on constructing classifying spaces was obtained by Rim\'anyi and Sz\H{u}cs~\cite{RS} who constructed a space $B_{\J}$ for each finite family $\J$ of {\it simple stable} singularities of maps of negative dimension. Remarkably, it lead Rim\'anyi to the method of Restriction Equations \cite{Ri} which prompted a series of explicit computations by Rim\'anyi, Kazarian, Feh\'er,  and others (e.g., see \cite{Ri}, \cite{FR}, \cite{Ka3} and references there). 

The Rim\'anyi-Sz\H{u}cs construction uses a technical J\"anich theorem~\cite{Ja},~\cite{Wa}, which a priori has no analogue  in the positive dimension case. However, when exists, the Rim\'anyi-Sz\H{u}cs construction $B_{\J}$ consists of classifying spaces $\BDiff F_{\alpha}$ of diffeomorphism groups of fibers $F_{\alpha}$ of $\J$-maps, i.e., as a set $B_{\J}$ is given by the disjoint union 
\[
   B_{\J}\cong \bigsqcup \BDiff F_{\alpha}.
\]
In the case $d> 1$, some of these strata are complicated. For example, the space $B_{\J}$ always contains a stratum $\BDiff F$ for each smooth manifold $F$ of dimension $d$; we recall that computing cohomology groups $H^*(\BDiff F; \Z)$ is a challenging (open) problem  already for a surface $F$ of a high genus. 

In contrast, thanks to an extremely helpful observation of Kazarian (e.g., see \cite{Ka2}), the proposed spectrum $\B_{\J}$ has a relatively simple structure both in the positive and negative dimension cases. It consists of Thom spaces of vector bundles over $S_t(\J)$, while the space $\mathbf{S}(\J)=\lim S_t(\J)$ has a natural stratification,
\[
     \mathbf{S}(\J)\cong \bigsqcup_{\tau\in \J} \BDiff {\tau},  
\]
where $\Diff \tau$ stands for the symmetry group of $\tau$. The spaces $\BDiff \tau$ are relatively simple. For example, in the case where $\J$ is the empty set of singularities of maps of dimension $d>0$, the space $\mathbf{S}(\J)$ coincides with $BO_d$ (compare with the corresponding monstrosity $\sqcup \BDiff F_{\beta}$, where the disjoint union ranges over closed manifolds $F_{\beta}$ of dimension $d$). The simplicity of $\B_{\J}$ makes it possible to carry out explicit computations not only in the case of $d<0$, but also in the case $d\ge 0$  \cite{Sz1}, \cite{Sad4}.     

\begin{remark}\label{r:1.10} Note that since both $\Omega^{\infty}\B_{\J}$ and $B_{\J}$ are classifying spaces, it follows that for each closed manifold $W$, there is a canonical isomorphism of sets
\[
   [W, \Omega^{\infty}\B_{\J}]\cong [W, B_{\J}]
\]
of homotopy classes of maps. In a paper in progress we study diffeomorphism groups of manifolds by exploring the relationship between $\Omega^{\infty}\B_{\J}$ and $B_{\J}$ from geometric point of view.
\end{remark}
 
\begin{remark} The b-principle is not a direct consequence of the h-principle. In terms of jet spaces, the h-principle for a differential relation $\mathcal R$ asserts that if $TV$ and $TW$ denote the tangent bundles of smooth closed manifolds $V$ and $W$ respectively, then the existence of a {\it formal solution} $V\to J^k(TV, TW)$ of $\mathcal R$ covering a map $f\co V\to W$ implies the existence of a homotopy of $f$ to a genuine solution of $\mathcal R$ (for definitions, see sections \ref{s:SDR}, \ref{s:JVB} and \ref{h-principle}); while the corresponding version of the  b-principle implies that for a sufficiently big integer $l\in \Z$ and stabilized tangent bundles $TV\oplus l\varepsilon$ of $V$ and $TW\oplus l\varepsilon$ of $W$, the existence of a {\it stable formal solution} $V\to J^k(TV\oplus l\varepsilon, TW\oplus l\varepsilon)$ of $\mathcal R$ covering the map $f\co V\to W$ implies the existence of a cobordism of $f$ to a genuine solution of $\mathcal R$ (for definitions, see sections \ref{s5}, \ref{s:JVB} and \ref{ssssss}). Here $\varepsilon$ stands for the trivial vector bundle over any topological space. 
\end{remark} 

\begin{remark} The proof of the Main Theorem in the present paper essentially utilizes the h-principle. The early versions of the h-principle appeared in the works of Nash, Smale, and later Hirsch, Poenaru, Phillips, Feit and others (e.g., see \cite{Gr}, \cite{EM} and references there). Its general version first appeared in the papers of Gromov \cite{Gr71} and Eliashberg and Gromov  \cite{EGr}, and in the present form was first formulated in the foundational book by Gromov \cite{Gr}. The theory has been extensively developed and several general powerful methods have been discovered (e.g., the Gromov's methods of convex integration, continuous sheaves \cite{Gr} and its new version, namely, the method of Eliashberg and Mishachev of holonomic approximations~\cite{EM}). We will use the  h-principle for so-called open $\mathcal{K}$-invariant differential relations which is essentially due to Phillips~\cite{Ph}, Eliashberg~\cite{El1}, \cite{El2}, du Plessis~\cite{Pl} and Ando~\cite{An6, An5, An4}. 
\end{remark}

\begin{remark}
It is plausible that a weak version of the h-principle is sufficient for proving the Main Theorem. Indeed, for finite families $\J$ of simple singularities in the negative dimension case~\cite{Sz1} one only needs the Rourke-Sanderson argument~\cite{RoS}; while in the case of maps of positive dimension a theorem similar to the Main Theorem can be established by using only the Gromov h-principle over open manifolds~\cite{Sad3}.    
\end{remark} 


\subsection*{Related constructions} The construction of $\Omega^{\infty}\B_{\J}$ is related to the construction by Eliashberg~\cite{El} of classifying spaces for Lagrangian and Legendrian immersions. In fact, our intermediate Theorem~\ref{t:iso} is a generalization of the Eliashberg theorem. 

According to a startling observation of Kazarian, the infinite loop space $\Omega^{\infty+d}\mathbf{MO}$ of the Thom spectrum of the unoriented cobordism groups contains a copy $K_{\J}$ of each space $B_{\J}$ \cite{Ka1}, \cite{Ka2} (see also Saeki-Yamamoto~\cite{SY}). Each space $K_{\J}$ comes with a natural stratification resembling that in the Rim\'anyi-Sz\H{u}cs construction. In the case $d<0$, the advantage of the classifying space $\Omega^{\infty}\B_{\J}$ over $K_{\J}$ is relatively limited, but in the case $d\ge 0$, the strata of $\Omega^{\infty}\B_{\J}$ are essentially simpler than those of $K_{\J}$ (see the discussion above).    

In the case of immersions, our construction should be compared to that of Wells~\cite{We}. 
In the case $d<0$ the Main Theorem is also proved by Sz\H{u}cs \cite{Sz1} by a different line of reasoning. In \cite{An4}, Ando proposes a homotopy theoretic counterpart of $B(W; \J)$ similar to but essentially different from ours and independently proves a theorem similar to but essentially different from our Theorem~\ref{t:iso}. It implies a version of the Main Theorem but only in the case $d<0$. Ando also suggests an alternative way of deriving the Main Theorem from our Theorem~\ref{t:iso} (cf. our Theorem~\ref{th:2} which originally appeared in~\cite{Sad5}). Again, the main advantage of our approach is that it applies not only in the case $d<0$, but also in the case $d\ge 0$ of our primary interest.

For $\J=\emptyset$ and $d=0$, the space $\Omega^{\infty}\B_{\J}$ can be identified with a path component of the infinite loop space $\Omega^{\infty}S^{\infty}$ and appears in the Kahn-Priddy theorem~\cite{Fu}. 
In the case where $\J=\emptyset$ and $d=2$, the space $\Omega^{\infty}\B_{\J}$ coincides with the space $\Omega^{\infty}\mathbf{hV}$ in the proof of the Mumford conjecture by Madsen and Weiss~\cite{MW}. In the case where $\J=\emptyset$ and $d>1$,  the space $\Omega^{\infty-1}\B_{\J}$ is weakly homotopy equivalent to the classifying space $BC_d$ of the cobordism category of manifolds of dimension $d$ of Galatius, Madsen, Tillman and Weiss \cite{GMTW}. In general $\Omega^{\infty-1}\B_{\J}$ is weakly homotopy equivalent to the classifying space $BC_{\J}$ of the corresponding singular cobordism category, provided $d>1$ \cite{Sad3}.

\section{Bordism principle}\label{s:2}

We may say that the homotopy principle for a set valued
functor $F$ asserts that $F$ can be induced from a homotopy functor with domain in a homotopy category. Likewise we
say that the {\it bordism principle} for a functor $F$ with values in
the category $\mathbf{AG}$ of abelian groups asserts that $F$ can be induced from a (co)homology functor.

We examine the b-principle in the case of $\mathbf{AG}$ valued functors $\EB$ of $\J$-bordism groups and $\EC$ of $\J$-cobordism groups.  The functors $\EB$ and $\EC$ are sometimes confused in the literature; both are defined so that for each closed manifold $W$ of a fixed dimension $n$, there is a canonical isomorphism
\[
     \EB(W)\cong \EC(W)\cong B(W; \J).
\]  
There is, however, an essential difference between $\EB$ and $\EC$. The former is covariant, while the latter is contravariant (see sections~\ref{s:5}, \ref{s:6}). 


We will show that under the conditions of the Main Theorem the covariant functor $\EB$ extends over the category $\Top^2\downarrow BO_n$ whose object is a vector $n$-bundle over a pair of topological spaces, and whose morphism is a fiberwise isomorphism of vector bundles. As we will see, the extended functor satisfies an analogue of the Eilenberg-Steenrod axioms. In other words, for the functor $\EB$ of $\J$-bordism groups we will establish a weak version of the $b$-principle (see Definition~\ref{def:2}). 

On the other hand, we will show that in contrast to the functor of $\J$-bordism groups, under the conditions of the Main Theorem, the functor of $\J$-cobordism groups does satisfy the b-principle, i.e., it can be induced from a genuine cohomology functor.    

We will need the following definitions. 

For a topological space $B$, let $\Top^2\downarrow B$ denote the
category with objects $(X,A;\varphi)$, each given by a pair of topological
spaces $(X,A)$ together with a homotopy class of a map $\varphi\co  X\to B$; and with
morphisms $(X_1,A_1; \varphi_1)\to (X_2, A_2; \varphi_2)$, each given by the homotopy class of a continuous map 
$f\co  (X_1, A_1)\to (X_2, A_2)$ such that $\varphi_1=\varphi_2\circ f$. There is a functor
$r\co \Top^2\downarrow B\to \Top^2\downarrow B$ defined by $r(X,A;
\varphi)=(A,\emptyset;\varphi|A)$.

\begin{definition}\label{def:1} A {\it homology theory $h_*$ on the category
$\Top^2\downarrow B$} is a sequence of functors $h_n\co 
\Top^2\downarrow B\to \mathbf{AG}$, with $n\in \Z$, and
natural transformations $\partial_n\co  h_n\to h_{n-1}\circ r$
that satisfy the Exactness axiom and the Excision axiom (e.g., see~\cite[Chapter 7]{Sw}).
\end{definition}

In a similar fashion one may give a definition of a cohomology theory $h^*$ on
the category $\Top^2\downarrow B$.

The category $\Top^2\downarrow B$ contains a subcategory $\Top\downarrow B$ of topological spaces $(X; \varphi)=(X,\emptyset; \varphi)$ over $B$ and homotopy classes of continuous maps over $B$.

\begin{definition}\label{def:2} Let $F$ be a contravariant (respectively covariant) functor
from a category $\mathcal{C}$ into the category of abelian groups
$\mathbf{AG}$. We say that the functor $F$ satisfies the {\it weak bordism principle} with respect to $\Top^2\downarrow B$ if there
is a cohomology (respectively homology) theory $h^*$  (respectively
$h_*$) on $\Top^2\downarrow B$ and a functor $\tau\co  \mathcal{C}\to \mathcal \Top\downarrow B$ such
that the functor $F$ is naturally equivalent to the functor
$h^n\circ \tau$ (respectively $h_n\circ \tau$) for some $n\in \Z$. We say that $F$ satisfies the {\it bordism principle} if $F$ satisfies the weak bordism principle in the case where $B$ is a point.  
\end{definition}

For example, the de Rham cohomology functor defined on the category of smooth manifolds and smooth maps satisfies the bordism
principle. Its homotopy analogue is the singular cohomology functor
defined on the category $\Top^2\downarrow \{pt\}$ of pairs of
topological spaces over a point $\{pt\}$.

\begin{remark} Since reduced and unreduced (co)homology are essentially
equivalent (e.g., see \cite{Hat}), Definition~\ref{def:1.1} given in the Introduction is essentially the same as Definition~\ref{def:2} of b-principle for contravariant functors.
\end{remark}

\section{Differential relations}\label{s:rel}\label{s:SDR}

Let $X\to V$ be a smooth fiber bundle. We
say that two smooth sections $f_1$ and $f_2$, defined in a
neighborhood of a point $v\in V$, have a {\it contact} of order
$\ge k$ at $v$ if the values and the partial derivatives of order
$\le k$ of $f_1$ and $f_2$ at $v$ are the same. The equivalence
class of local sections that have a contact of order $k$ at $v$
with a local section $f$ is called the {\it $k$-jet $[f]_v^k$}.
The set of $k$-jets of local sections of $X$ forms the total space $J^k(X)$, called the {\it $k$-jet space}, of a smooth fiber bundle over $X$ with projection $\pi_X$ sending a $k$-jet
$[f]_v^k$ onto $f(v)$. The composition of $\pi_X$ with the
projection $X\to V$ turns $J^k(X)$ into the total space of a smooth
fiber bundle over $V$.

\begin{definition}  A {\it differential relation $\mathcal R$} of order $k$ over a smooth fiber bundle $X\to V$ is a subset of the $k$-jet space
$J^k(X)$.
\end{definition}

Every smooth section $f$ of the bundle $X\to V$ is covered by a
smooth map
\[
  j^kf\co  V\to J^k(X),
\]
\[
  j^kf\co  v \mapsto [f]_v^k,
\]
called the {\it $k$-jet extension of $f$}.

\begin{definition} A {\it solution} of a
differential relation $\mathcal R\subset J^k(X)$ is a section $f$ with
image $j^kf(V)$ in $\mathcal R$.
\end{definition}

In other words, a differential relation over a smooth fiber bundle is a relation on the derivatives of smooth sections of the bundle. 

If $X\to V$ is a smooth fiber bundle over a manifold of dimension $m$
with fiber of dimension $n$, then the fiber of the $k$-jet bundle
$J^k(X)\to X$ is isomorphic to the space $J^k(\R^m,\R^n)$ of
$k$-jets of germs
\begin{equation}\label{eq:germ}
f\co  (\R^m,0) \longrightarrow (\R^n,0).
\end{equation}
Occasionally we will identify the space $J^k(\R^m, \R^n)$ with the space of polynomial maps $(\R^m, 0)\to (\R^n, 0)$ of degree $\leq k$.

The space $J^k(\R^m, \R^n)$ admits a smooth action of the so-called {\it
$k$-contact group} $\K=\K(k,m,n)$. By definition, the group $\K$ is the subgroup of $k$-jets of diffeomorphism germs
\[
    (\R^m\times \R^n, 0) \longrightarrow (\R^m\times \R^n, 0)
\]
consisting of those elements which take the horizontal slice $\R^m\times \{0\}$ onto $\R^m\times \{0\}$ and each vertical slice $\{x\}\times \R^n$, with $x\in \R^m$, into a vertical slice
$\{y\}\times \R^n$, with $y\in \R^m$. The action of $\K$ on $J^k(\R^m, \R^n)$ is determined by the action of $\K$ on
the graphs $\{(x, f(x))\in \R^m\times \R^n\}$ of germs of the form
(\ref{eq:germ}).


\begin{definition} A {\it differential relation basis $R$} of order
$k$ is an arbitrary subset of the space $J^k(\R^m,\R^n)$. We say that a basis {\it has sufficiently many symmetries} or {\it $\K$-invariant} if it is invariant with respect to the action of $\K$.
\end{definition}

A basis invariant with respect to the action of $\K(k,m,n)$, leads to a
differential relation over every trivial fiber bundle $\pi_V\co 
X=V\times W\to V$ with base of dimension $m$ and fiber of
dimension $n$. Indeed, over the total space $X$, there is a
canonically defined, principal $\K$-bundle $P^k_X(\K)\to X$ of $k$-jets of diffeomorphism germs 
\begin{equation}\label{eq:3}
   \alpha\co  (\R^m\times \R^n, 0) \longrightarrow (X, x), \qquad x=(v,w)\in X,
\end{equation}
such that the restriction $\alpha|(\R^m\times \{0\})$ is a
diffeomorphism onto a neighborhood of $(v,w)$ in $V\times \{w\}$, and for each point $y\in \R^m$, the image of the vertical slice
$(\{y\}\times \R^n)$ under the map $\pi_V\circ \alpha$ is a single
point in $V$. The projection $P^k_X(\K)\to X$ is given by taking the $k$-jet represented by a germ (\ref{eq:3}) onto $x$. Let us observe that there is a canonical isomorphism of fiber bundles
\[
    P^k_X(\K)\times_{\K} J^k(\R^m, \R^n) \longrightarrow J^k(X)
\]
over $X$, which, in terms of representatives, takes a pair $(\alpha, f)$ of germs (\ref{eq:3}) and (\ref{eq:germ}) onto a local section with graph 
\[
        \alpha\circ ([\pi_V\circ \alpha|({\R^m\times\{0\}})]^{-1},\ f \circ [\pi_V\circ \alpha|({\R^m\times\{0\}})]^{-1}).
\]
If, now, $R\subset J^k(\R^m, \R^n)$ is a $\K$-invariant basis, then the set $P^k_X(\K)\times_{\K} R$ is a well-defined relation over $X$.

%

\begin{remark} The sections of a trivial bundle $X=V\times W\to V$ are in bijective correspondence with the maps $V\to W$. In this case, according
to the terminology of singularity theory, the complement $\Sigma =
J^k(X) \setminus \mathcal R$ to a differential
relation $\mathcal R$ is a {\it singularity set}, and a solution to
$\mathcal R$ is a {\it map without $\Sigma$-singularities}. Thus a choice of a $\K$-invariant differential relation basis containing regular jets corresponds to a choice of a (not necessarily open) $\K$-invariant set of prescribed singularity types $\J$ (see section~\ref{s:1}). Namely, the set $\J$ corresponding to $\mathcal R$ consists of non-regular singularity types of germs whose jets are in $\mathcal R$. We note that if the set of singularity types $\J$ is open, then the corresponding differential relation basis $R$ is open in the topological space $J^k(\R^m, \R^n)$. 
\end{remark}

\section{Sequence of differential relation bases}\label{s5}

We have seen that a base $R\subset J^k(\R^m, \R^n)$ with sufficiently many symmetries determines a differential relation on smooth maps $V\to W$ for any pair of manifolds $V$ and $W$ of dimensions $m$ and $n$ respectively.  However, to define a bordism group of solutions we also need to consider differential relations on maps of manifolds of higher dimensions, which can be done by introducing a {\it sequence of differential relation bases}.

For non-negative integers $s$, $t$, with $s<t$, and an integer $q$,  there is an embedding 
\[
ext^s_t\co  J^k(\R^{s}, \R^{q+s}) \longrightarrow J^k(\R^{t},
\R^{q+t})
\]
of $k$-jet spaces that takes the $k$-jet of a germ
\[
f\co (\R^{s},0)\to (\R^{q+s},0)
\]
onto the $k$-jet of the suspension germ
\[
f\times \id_{\R^{t-s}}\co  (\R^{t},0)\to (\R^{q+t},0),
\]
where $\id_{\R^{t-s}}$ stands for the identity map of the space $\R^{t-s}$. We note that the space $J^k(\R^s, \R^{q+s})$ is empty if $q+s< 0$.

\begin{definition} A {\it sequence $\mathbf{R}=\mathbf{R}(q)$ of differential relation bases}, or a {\it stable differential relation basis} is a set of bases $R_s(q)\subset J^k(\R^{s}, \R^{q+s})$, one for each $s\ge 0$, such that for non-negative integers $s,t$, with $s<t$, the map $ext^s_t$ takes $R_s=R_s(q)$ into $R_t$ and the complement of $R_s$ in $J^k(\R^{s}, \R^{q+s})$
into the complement of $R_t$ in $J^k(\R^{t}, \R^{q+t})$. We say that a sequence $\mathbf{R}=\{R_s\}$ of bases  is {\it open} if each basis $R_s$ is open in $J^k(\R^s, \R^{q+s})$.
\end{definition}

\begin{remark} For a map $f\co  V\to W$ of a manifold of dimension $t$ to a manifold of dimension $t+q$, the number $q$ is called the {\it codimension} of $f$. Thus, if $f$ is of dimension $d$, then $q=-d$. We tend to use the term ``codimension" if $q\ge 0$, e.g., in the case of immersions. On the other hand, in the case $q<0$ we find the term ``dimension" more appropriate; for example, the phrase ``a map $f$ of dimension 2" means that a non-empty regular fiber of $f$ is a manifold of dimension $2$, while the phrase ``a map of codimension $-2$" has no obvious geometric meaning.   
\end{remark}

We may identify a sequence $\mathbf{R}(q)$ of bases with a subset of the space
\[
   \mathbf{J}=\mathbf{J}(q)\co =\colim_{s\to\infty}J^k(\R^{s},\R^{q+s}),
\]
where the colimit is taken with respect to the maps
$ext^s_{s+1}$, which we regard as inclusions. Indeed, if $\mathbf{R}(q)$ is a
subset of $\mathbf{J}(q)$, then the set $\{R_s\}$ of spaces
$R_s=\mathbf{R}\cap J^k(\R^{s},\R^{q+s})$ is a sequence of
differential relation bases. Conversely, if $\{R_s(q)\}$ is a sequence of
differential relation bases, then the colimit 
\[
   \mathbf{R}=\mathbf{R}(q)\co =\colim_{s\to\infty}R_s(q)
\]
taken with respect to the maps $ext^s_{s+1}$ is a subset
of $\mathbf{J}(q)$. 

Similarly there are natural inclusions $\K(k,s,q+s)\to \K(k,t,q+t)$, $s<t$, of groups, which allow us to define a stable $k$-contact group 
\[
   \mathbf{K}=\mathbf{K}(q)\co =\colim_{s\to\infty}\K(k, s, q+s).
\]
The actions of the groups $\K(k,s,q+s)$ on $J^k(\R^{s}, \R^{q+s})$ induce an action of the group $\mathbf{K}$ on $\mathbf{J}$. It follows that a stable basis $\mathbf{R}\subset \mathbf{J}$ is {\it $\mathbf{K}$-invariant} if and only if each relation basis $R_s$, with $s\ge 0$, is invariant with respect to the action of $\K(k, s, q+s)$. In this case we will also say that $\mathbf{R}$ has {\it sufficiently many symmetries}. 

A sequence of bases $\{R_s\}$ is a {\it suspension} of a (not necessarily $\mathcal{K}$-invariant) differential relation basis $R$ if $R_t\supset R$ for an appropriate $t$. A suspension is a {\it $\mathbf{K}$-suspension} if
the sequence $\{R_s\}$ is $\mathbf{K}$-invariant. We observe that each basis $R\subset J^k(\R^s, \R^{q+s})$ has the minimal $\mathbf{K}$-suspension, defined as the minimal $\mathbf{K}$-invariant subset of $\mathbf{J}$ that contains $R$; and the maximal $\mathbf{K}$ suspension, defined as the maximal $\mathbf{K}$-invariant subset of $\mathbf{J}$ that does not intersect the set $J^k(\R^s,\R^{q+s})\setminus R$. 

\begin{remark} If $R$ is a differential relation basis, then its minimal $\mathbf{K}$-suspen\-sion exists even if $R$ is not $\K$-invariant. On the other hand the maximal $\mathbf{K}$-suspension of $R$ exists only if $R$ is $\K$-invariant. Indeed, suppose that $\{R_s\}$ is the maximal suspension of $R$. Since it is a suspension, there is an inclusion $R_t\supset R$ for some $t$. On the other hand, since the suspension is maximal, there is an inclusion $R_t\subset R$. Consequently, $R_t=R$.  
\end{remark}

\begin{example}[Immersions] If $R=R_0, R_1,...,$ stand for
the bases corresponding to immersions in
respectively $J^k(\R^0,\R^q)$, $J^k(\R^{1},\R^{q+1}), ...$, then
$\{R_s\}$ is the minimal suspension of $R\subset J^k(\R^0, \R^q)$.
\end{example}

\begin{remark} Due to Mather it is well-known that if $\{R_s\}$ is the minimal suspension of a $\mathcal{K}$-invariant basis $R\subset R_t$, then $R_t=R$. On the other hand, we may define the minimal $G$-suspension $\{R_s\}$ of $R$ for any subgroup $G<\mathbf{K}$, e.g., for the group $\mathcal A$ of $k$-jets of right-left coordinate changes. In this case an $\mathcal A$-invariant basis $R\subset R_t$ may not coincide with the $t$-th space $R_t$ of the minimal $\mathcal A$-suspension of $R$. 
\end{remark}

\section{Bordism groups of solutions}\label{s:5}

A sequence $\mathbf{R}$ of bases $R_s\subset J^k(\R^{s}, \R^{q+s})$, $s\ge 0$, with sufficiently many symmetries determines a differential relation $\mathcal{R}$ for mappings from any manifold
$V$ of dimension $s$ into any manifold $W$ of dimension $q+s$.
To simplify formulation of statements, we will say that a map satisfying
$\mathcal{R}$ is a {\it solution of $\mathbf{R}$}, or, simply, an {\it  $\mathbf{R}$-map}.

By definition, two $\mathbf{R}$-maps $f_i\co  V_i\to W_i$, with $i=1,2$, of closed manifolds are {\it right-left bordant} if there are
\begin{itemize}
\item a compact manifold $V$ with $\partial V= V_1\sqcup V_2$,
\item a compact manifold $W$ with $\partial W=W_1\sqcup W_2$, and
\item an $\mathbf{R}$-map $f\co V\to W$ such that $f(V_i) \subset
W_i$ for $i = 1, 2$, and the restriction
of $f$ to collar neighborhoods of $V_1$ and $V_2$ in $V$ can be identified with the disjoint union of suspensions of $f_1$ and $f_2$.  
\end{itemize}
The right-left bordism classes of  $\mathbf{R}$-maps constitute a group with addition defined in terms of representatives by taking the disjoint union of maps. Each element of the right-left bordism group is of order $2$. 

We say that two $\mathbf{R}$-maps $f_i\co V_i\to W$, with $i=1,2$, of closed manifolds $V_i$ into a manifold $W$ are {\it (right) bordant} if there are
\begin{itemize}
\item a compact manifold $V$ with $\partial V= V_1\sqcup V_2$, and
\item an $\mathbf{R}$-map $f\co V\to W\times [1,2]$ such that $f(V_i)\subset W\times \{i\}$ for $i=1,2$, and the restriction of $f$ to collar neighborhoods of $V_1$ and $V_2$ in $V$ can be identified with the disjoint union of suspensions of $f_1$ and $f_2$.  
\end{itemize}

\begin{remark} Note that in contrast to the definition of right-left bordisms, in the definition of right bordisms the target manifold $W$ is not assumed to be closed. For our approach, in the definition of right bordisms it is essential to allow arbitrary (smooth) manifolds (without boundary), not only closed manifolds. Indeed, Theorem~\ref{t:iso} implies, for example, that the right bordism groups of $\mathbf{R}$-maps to non-closed manifolds $\R^n$, $n\ge 0$, correspond to terms of what in homology theory is called {\it the coefficient group}.    
\end{remark}

Finally, two $\mathbf{R}$-maps $f_i\co V_i\to W_i$, with $i=1,2$, of closed manifolds are {\it left bordant} if $V_1=V_2$ and there
is a right-left bordism $(V,W,f)$ with $V=V_1\times [0,1]$.

Taking the disjoint union of maps leads to a structure of a semigroup on the set of right bordism classes of solutions. We note that the semigroup of bordism classes of solutions may not be a group. 

\begin{remark} We adopt the convention that the empty set $\emptyset$ is a manifold of an arbitrary dimension. In particular, for each manifold $W$ and a sequence $\mathbf{R}$ of relation bases, the map $\emptyset\to W$ is a solution of $\mathbf{R}$.   
\end{remark}

\begin{example} Let $\mathbf{R}$ be the sequence of open differential relation bases $R_s$, with $s\ge 0$, with sufficiently many symmetries corresponding to submersions of dimension $-q>0$. Then the element of the semigroup of bordism classes of solutions represented by any submersion $V\to W$ of a closed non-empty manifold has no inverse. Consequently, the semigroup of bordism classes of submersions is not a group.  
\end{example}

We define the {\it (right) bordism group} of
$\mathbf{R}$-solutions $\EB(W)$ as the group given by the Grothendieck construction applied to the semigroup $SB(W)$ of right bordism classes of solutions.
In fact, we will see that under the assumptions of the Main Theorem the semigroup of  right bordism classes of solutions is already a group (see Corollary~\ref{c:1}). Furthermore, one may follow the proof of Theorem~\ref{t:iso} to obtain a description of the inverse element (see \cite{Sz1} for the case $q>0$), but an explicit construction of a representative of the inverse of a given element in $SB(W)$ may not however  be a simple task.

Let $\mathbf{Diff}_n$ be the category of smooth manifolds of dimension $n$ without boundary and equidimensional embeddings. We have defined the correspondence
\[
       \EB\colon \mathop{\mathrm{Obj}}(\mathbf{Diff}_n)\longrightarrow \mathop{\mathrm{Obj}}(\mathbf{AG}),
\]
\[
       \EB\colon W \mapsto \EB(W),
\]
where $\mathop{\mathrm{Obj}}(C)$ stands for the collection of objects of a category $C$. The collection of morphisms of a category $C$ will be denoted by $\mathop{\mathrm{Mor}}(C)$.
 
Given an embedding $i\co  W_1\to W_2$ in $\mathop{\mathrm{Mor}}(\mathbf{Diff}_n)$, and an $\mathbf{R}$-map $f\co  V_1\to W_1$ representing an element $[f]$ in $\EB(W_1)$, the composition $i\circ f\co  V_1\to W_2$ is an $\mathbf{R}$-map representing an element $[i\circ f]$ in $\EB(W_2)$. Furthermore, the element $[i\circ f]$ depends only on the class $[f]$, not on its representative $f$. Hence, each embedding $i\co  W_1\to W_2$ in $\mathop{\mathrm{Obj}}(\mathbf{Diff}_n)$ gives rise to a correspondence 
\[ 
        \EB(i)\colon \EB(W_1) \longrightarrow \EB(W_2),
\]
\[
         \EB(i)\colon [f]-[g]\mapsto [i\circ f]-[i\circ g],
\] 
which is easily seen to be a group homomorphism. 
In fact the correspondence $\EB$ defines an $\mathbf{AG}$ valued covariant functor on the category $\mathbf{Diff}_n$, called the {\it functor of the $\mathbf{R}$-bordism group}. We will continue to denote the functor by the symbol $\EB$. 

\begin{remark} We warn the reader that some authors use the term ``cobordism groups" to refer to the values of the functor $\EB$ and define ``bordism groups" as cobordism groups of Euclidean spaces. Other authors use the term ``cobordism groups" to refer both to groups related to the covariant functor $\EB$ and to groups related to the contravariant functor $\EC$ described in section~\ref{s:6}. 

We use the term ``bordism groups" in the case of covariant functors and reserve the term ``cobordism groups" for contravariant functors (see section~\ref{s:6}). Our choice of terminology is compatible with that in \cite{Sw}, \cite{Ru}, \cite{St}, \cite{Koch} and other textbooks in algebraic topology.    
\end{remark}

\section{Cobordism groups of solutions}\label{s:6}

Two proper $\mathbf{R}$-maps $f_i\co V_i\to W$, with $i=1,2$, of manifolds are {\it (right) cobordant} if there are
\begin{itemize}
\item a manifold $V$ with $\partial V= V_1\sqcup V_2$, and
\item a proper $\mathbf{R}$-map $f\co V\to W\times [1,2]$ such that $f(V_i)\subset W\times\{i\}$ for $i=1,2$, and the restriction of $f$ to collar neighborhoods of $V_1$ and $V_2$ in $V$ can be identified with the disjoint union of suspensions of $f_1$ and $f_2$.  
\end{itemize}
Taking the disjoint union of maps leads to a structure of a semigroup on the set of  right cobordism classes of solutions. 

We recall that $\mathbf{Diff}_n$ denotes the category of smooth manifolds of dimension $n$ without boundary and equidimensional embeddings. 

The {\it (right) cobordism group} of
$\mathbf{R}$-solutions $\EC(W)$ is defined to be the group given by the Grothendieck construction applied to the semigroup of right cobordism classes of solutions. In other words, $\EC(\cdot)$ is a correspondence
\[
       \EC\colon \mathop{\mathrm{Obj}}(\mathbf{Diff}_n)\longrightarrow \mathop{\mathrm{Obj}}(\mathbf{AG}),
\]
\[
       \EC\colon W \mapsto \EC(W).
\]
Given an embedding $i\co  W_1\to W_2$ in $\mathop{\mathrm{Mor}}(\mathbf{Diff}_n)$, and an $\mathbf{R}$-map $f\co  V_2\to W_2$ representing an element $[f]$ in $\EC(W_2)$, the pullback 
\[
i^*f=i^{-1}\circ f\co  V_1=f^{-1}(i(W_1))\longrightarrow W_1
\] 
is an $\mathbf{R}$-map representing an element $[i^*f]$ in $\EC(W_1)$. The element $[i^*f]$ depends only on the class $[f]$, and therefore, each embedding $i$ in $\mathop{\mathrm{Mor}}(\mathbf{Diff}_n)$ determines a correspondence
\[ 
        \EC(i)\colon \EC(W_2) \longrightarrow \EC(W_1),
\]
\[
         \EC(i)\colon [f]-[g]\mapsto [i^*f]-[i^*g],
\] 
which is a group homomorphism. The correspondence $\EC$ defines an $\mathbf{AG}$ valued contravariant functor on $\mathbf{Diff}_n$, called the {\it functor of the $\mathbf{R}$-cobordism group}. We will continue to denote this functor by $\EC$. 


\begin{remark} As has been mentioned in section~\ref{s:2}, it follows from the definitions that if $W$ is a closed manifold, then there is an isomorphism 
\[
             \EB(W)\stackrel{\cong}\longrightarrow \EC(W).
\]
In fact, in terms of representatives the isomorphism is given by the identity correspondence
\begin{equation}\label{eq:6.1}
                 [f]\mapsto [f].
\end{equation}
On the other hand, if $W$ is not closed, then the correspondence (\ref{eq:6.1}) is not an isomorphism. For example, for any sequence $\mathbf{R}(q)$, with $-q>0$, the inclusion $f\co  \R^0\to \R^{-q}$ represents a non-trivial element in $\EB(\R^{-q})$ and the trivial element in $\EC(\R^{-q})$.  
\end{remark}

\section{Formal jet spaces of vector bundles}\label{s:JVB}
Given two vector bundles $\eta$ and $\gamma$ of dimensions $m$ and $n$ respectively over a topological
space $B$, let $P^k_{\eta,\gamma}(\Symm)$ denote the principal
$\Symm$-bundle over $B$ whose total  space consists of $k$-jets of diffeomorphism germs 
\begin{equation}\label{eq:2}
   f\co (\R^m\times \R^n, 0) \longrightarrow (\eta|x\oplus\gamma|x,
   0), \qquad x\in B,
\end{equation}
where $\eta|x$ and $\gamma|x$ are restrictions to $x$, such that $f(\R^m\times \{0\})=\eta|x\oplus 0$ and for each $y\in \R^m$, the projection
$\eta \oplus \gamma\to \eta$ maps the image $f(\{y\}\times \R^n)$
onto a point. The projection $P^k_{\eta,\gamma}(\Symm)\to B$ is defined by sending the germ (\ref{eq:2}) onto $x\in B$.

\begin{definition}\label{def:diff} A {\it formal differential relation $\mathcal R$} over $(\eta,\gamma)$ is an arbitrary subset of the total space of the {\it formal $k$-jet bundle} 
\[
   J^k(\eta,\gamma)=P^k_{\eta, \gamma}(\Symm)\times_{\Symm} J^k(\R^m,\R^n) \longrightarrow B.
\]
\end{definition}

Given a $\K$-invariant basis $R\subset J^k(\R^m,\R^n)$, there is
an associated differential relation
\[
   \mathcal R=\mathcal R(\eta,\gamma)=P^k_{\eta, \gamma}(\K)\times_{\K} R
\]
 in $J^k(\eta,\gamma)$. 

In fact, for a smooth trivial bundle $\pi_V\co  X=V\times W\to V$ where the base and fiber are Riemannian manifolds, the $k$-jet bundle $J^k(X)$ is canonically isomorphic to a certain formal $k$-jet bundle. In what follows, the tangent bundle of a manifold $M$ is denoted by $TM$ and the tangent plane of $M$ at $x$ is denoted by $T_xM$. The Riemannian metric on $X$ is given by the product of Riemannian metrics on $V$ and $W$.

\begin{lemma}\label{l:equivalence} Let $\eta$ be the pullback of $TV$ to the bundle over $X$ with respect to $\pi_V$ and $\gamma$ the subbundle of $TX$ that consists of vectors
tangent to the fibers of the bundle $X\to V$. Then there is a
canonical isomorphism $J^k(\eta, \gamma)\to J^k(X)$.
\end{lemma}
\begin{proof} Let $\Exp_x\co  U_x\to V_x$ denote the exponential
diffeomorphism of a neighborhood of the origin in the tangent
plane $T_xX$ onto a neighborhood of $x$ in $X$.

Given a map $\alpha\co  (\R^m\times \R^n, 0)\to (T_xX, 0)$
representing a point in $P^k_{\eta, \gamma}(\K)$, the composition
$\Exp_x\circ \alpha$ represents an element in $P^k_X(\K)$. Hence, the product Riemannian metric on $X$ determines
a map $P^k_{\eta, \gamma}(\K)\to P^k_X(\K)$ which is easily seen to
be an isomorphism of principal $\K$-bundles. The isomorphism of Lemma~\ref{l:equivalence} can be produced from the isomorphism of the principal $\K$-bundles by the Borel construction. 
\end{proof}

\section{Bordism group of formal solutions}\label{ssssss}

Let $\R^{\infty}$ denote the infinite dimensional vector space, defined  as the colimit of inclusions $\R^1\subset \R^2 \subset ...$, with base $(e_1,e_2,\dots)$. We regard the classifying space $BO_i$ as the space of
$i$-subspaces of $\R^{\infty}$ each of which is contained in a
subspace $\R^j$ for some $j<\infty$. Then the fiber of a
canonical vector bundle $EO_i\to BO_i$ over a plane $L\in BO_i$ can be interpreted as the space of vectors in $L$. There is a canonical inclusion $BO_i\to BO_{i+1}$ that takes a plane spanned by vectors $v_1,\cdots, v_i$ onto the plane spanned by $e_1,\theta(v_1),\cdots, \theta(v_i),$ where $\theta$ is the shift endomorphism of $\R^{\infty}$ given by $\theta(e_k)=e_{k+1}$. Similarly there is a canonical inclusion $EO_{i}\subset EO_{i+1}$ and splitting of $EO_{i+1}|BO_i$ into the sum $\varepsilon\oplus EO_i\cong EO_i\oplus \varepsilon$. Here and in what follows we use the
symbol $\varepsilon$ to denote the trivial $1$-dimensional vector
bundle over an arbitrary space.

Let $W$ be a manifold of dimension $n$. Then a $\mathbf{K}$-invariant sequence $\mathbf{R}=\{R_s(q)\}$ of bases in $J^k(\R^{s}, \R^{q+s})$, with $s\ge 0$ and $q+s\ge n$, leads to a sequence of formal differential relations $[R_W]_s$ in the formal $k$-jet bundles $[J^k_W]_s$ associated to the pair of bundles 
\[
(\pi_1^*EO_{s},\ \pi_2^*{TW\oplus t \varepsilon}) \quad \textrm{over}\quad BO_{s}\times W, 
\]
where $\pi_1$ and $\pi_2$ denote the projections of $BO_{s}\times W$ onto the first and second factors respectively, and $t=q+s-n$. 

For $s\ge 0$, there is a canonical inclusion $[J^k_W]_s\subset [J^k_W]_{s+1}$ that takes the $k$-jet $[f]^k_0\in [J^k_W]_s$ over $v\in BO_{s}\times W$ represented by the germ of a smooth map
\[
     f\co  (\pi_1|v)^*EO_{s} \longrightarrow (\pi_2|v)^*\{TW\oplus t\varepsilon\}
\]
onto the $k$-jet $[f\times \id_{\varepsilon}]^k_0\in [J^k_W]_{s+1}$ of the germ of the map
\[
     f\times \id_{\varepsilon} \co  (\pi_1|v)^*(EO_{s}\oplus \varepsilon) \longrightarrow (\pi_2|v)^*\{TW\oplus t\varepsilon\oplus\varepsilon\},
\]
where we make use of the identification of $EO_{s+1}|BO_s$ with $EO_s\oplus \varepsilon$.

Let us observe that, since $\mathbf{R}$ is a $\mathbf{K}$-invariant sequence of bases, the canonical inclusion $[J^k_W]_s\subset [J^k_W]_{s+1}$ takes $[R_W]_s$ into $[R_W]_{s+1}$. In view of these canonical inclusions, we will use the same symbol $\pi$ to denote  any of the projections
\[
     [J^k_W]_s\longrightarrow BO_{s}, \qquad s\ge 0,
\]
\[
     [R_W]_s\longrightarrow BO_{s}, \qquad s\ge 0,
\]
defined as the composition of the projection onto $BO_{s}\times W$ followed by the projection onto the first factor $BO_{s}$. 

\begin{definition}\label{def:3} Let $\mathbf{R}=\mathbf{R}(q)$ be a $\mathbf{K}$-invariant stable differential relation basis, $W$ a smooth manifold of dimension $n\ge q$, and $m=n-q$. Then, {\it the bordism group $\M$ of stable formal $\mathbf{R}$-maps} into $W$ is the $m$-th bordism group of maps of manifolds with an $(R_W,\pi)$-structure in the stable tangent bundle, where $(R_W, \pi)$ is the $(B,f)$-sequence (for a definition, see, for example, \cite{St}) given by the commutative diagrams 
\[
\begin{CD} [R_W]_s @>>> [R_W]_{s+1} \\
     @V\pi VV @V\pi VV \\
     BO_{s} @>>> BO_{s+1},
\end{CD}
\] 
indexed by $s\ge m$, with horizontal maps given by the canonical inclusions. 
\end{definition}

In particular an element of the bordism group of formal solutions is represented by a triple $(V, \alpha, \tilde\alpha)$ of a closed manifold $V$ of dimension $m$, a continuous map $\alpha\co  V\to [R_W]_s$ for some $s\ge m$, and a fiberwise isomorphism of vector bundles
\[
    \tilde\alpha\co  TV\oplus (s-m)\varepsilon \longrightarrow \pi^*EO_{s}
\] 
covering the map $\alpha$. In particular, the composition $\pi\circ \alpha$ is a map classifying the stable tangent bundle of $V$.

By the Pontrjagin-Thom construction, the $m$-th bordism group of manifolds with $(R_W, \pi)$-structure in the stable tangent bundle is isomorphic to the $m$-th homotopy group of the Thom spectrum, denoted by $\mathbf{MO\wedge_R} W$, associated with $(R_W, \pi)$ (see, for example, \cite{St}). 

Note that if $\mathbf{R}$ consists of all jets, i.e., if any map of codimension $q$ is an $\mathbf{R}$-map, then $\mathbf{MO\wedge_R} W\cong \mathbf{MO}\wedge W$ and Theorem~\ref{t:iso} below reduces to the classical Pontrjagin-Thom theorem. 

In section~\ref{s:11} we will show that $\mathbf{MO_R}(\cdot)$ can be extended to a covariant functor on the category $\Top^2 \downarrow BO_n$.

\section{H-principle}\label{h-principle}

Certainly, a necessary condition for the existence of a solution
of a differential relation $\mathcal R$ is the existence of a
section of the bundle $\pi_V\co J^k(X)\to V$ with image in $\mathcal
R$. The homotopy converse statement, which may not be true in general, is
referred to as the {\it homotopy principle in the existence level}, or, simply, {\it h-principle}.

\begin{definition}[H-principle] Every section $s\co  V\to
\mathcal R\subset J^k(X)$ is homotopic to the $k$-jet extension of
a section $V\to X$ by homotopy of sections $V\to \mathcal R$.
\end{definition}

In particular, if a differential relation satisfies the
h-principle, and there is a section $s\co V\to \mathcal R$, then
the differential relation has a solution.

By the Gromov theorem~\cite{Gr}, an open so-called $\Diff$-invariant differential relation over a fiber bundle $X\to V$ always satisfies the homotopy principle if $V$ is an open manifold. We refer the reader to \cite{Gr}, \cite{EM} and \cite{Sp} for further
examples of differential relations satisfying the h-principle. 

A section $V\to J^k(X)$ is {\it holonomic} if it extends a section
$V\to X$.

We will also need the relative version of the h-principle.

\begin{definition}[Relative h-principle] Suppose that a
section $s\co V\to \mathcal R$ restricted to a neighborhood of a
closed subset $V_0\subset V$ is holonomic. Then $s$ is homotopic to a
holonomic section through sections $V\to \mathcal R$ constant in
a neighborhood of $V_0$. 
\end{definition}

Under the conditions of the Main Theorem, the relative h-principle for open $\mathcal K$-invariant differential relations imposed on maps of manifolds always holds true. It is essentially due to Phillips~\cite{Ph}, Eliashberg~\cite{El1}, \cite{El2}, du Plessis~\cite{Pl} and Ando~\cite{An6, An5, An4}.

We say that a $\mathbf{K}$-invariant sequence $\mathbf{R}=\{R_s\}$ of differential relation bases $R_s\subset
J^k(\R^{s}, \R^{q+s})$, $s\ge 0$, {\it satisfies the (relative)
h-principle} if for every trivial bundle $X\to V$ over an $s$-dimensional
manifold with a $(q+s)$-dimensional fiber, the corresponding
differential relation in $J^k(X)$ satisfies the (relative)
h-principle.

\section{Destabilization argument}\label{s:iso}

Let $\mathbf{R}$ be a sequence of differential relation bases $R_s\subset J^k(\R^{s}, \R^{q+s})$, $s\ge 0$, with sufficiently many symmetries. Then, for any smooth manifold $W$ of dimension $n$, there is a semigroup homomorphism from 
the semigroup $SB(W)$ of right bordism classes of $\mathbf{R}$-maps into $W$ to the bordism group $\M$ of stable formal $\mathbf{R}$-maps into $W$. 

Indeed, suppose that $V$ is a manifold of dimension $m=n-q$ and $f\co V\to W$ is a smooth map representing an element in the semigroup $SB(W)$. Let $\eta$ and $\gamma$ be the vector bundles over $X=V\times W$ defined as the pullbacks of the tangent bundles of $V$ and $W$ with respect to the projections of $V\times W$ onto the first and second factors respectively. Then, in view of  Lemma~\ref{l:equivalence}, by choosing Riemannian metrics on $V$ and $W$, we may identify $J^k(\eta, \gamma)$ with the $k$-jet space $J^k(X)$ of the trivial bundle $X\to V$. In particular, we may say that the map $f$ gives rise to a section $j^kf$ of the fiber bundle $\pi_V^k\co  \mathcal{R}_m(\eta, \gamma)\to V$.

Let $\tau\co  TV\to EO_m$ be a bundle map classifying the tangent bundle of $V$. In several occasions we will implicitly use a well-known fact that the homotopy class of $\tau$ is unique. 


The bundle map $\tau$ gives rise to a fiberwise isomorphism of bundles 
\[
\begin{CD}
J^k(\eta, \gamma) @>>> [J_W^k]_m   \\
@V{\pi_V^k} VV @VV\pi V \\
V @>>> BO_m.
\end{CD}
\]
Let $\tilde\alpha\co  TV\to \pi^*EO_m$ be the canonical lift of $\tau$ covering the composition $\alpha$ of $j^kf$ and the fiberwise isomorphism $J^k(\eta, \gamma)\to [J_W^k]_m$. Then the triple $(V, \alpha, \tilde\alpha)$ represents an element in $\M$, which, as it is easy to see, does not depend on $\tau$ and the choice of Riemannian metrics on $V$ and $W$. Furthermore, if $f_i\co V_i\to W$, $i=1,2$, are two $\mathbf{R}$-maps representing the same element
in $SB(W)$, \ie, if there exists a bordism between $f_1$ and $f_2$
satisfying $\mathbf{R}$, then, by a similar argument, the triples corresponding to $f_1$ and $f_2$ determine the same element in $\M$. Thus, there is a well-defined map
\[
   \psi\co  SB(W)\longrightarrow \M,
\]
which is easily seen to be a semigroup homomorphism. 

%
%
%

\begin{theorem}\label{t:iso} The homomorphism $\psi$ is an isomorphism for any (not necessarily closed) manifold $W$ of dimension $n$. In particular, the bordism semigroup $SB(W)$ of $\mathbf{R}$-maps into $W$ is isomorphic to the group $\pi_m(\mathbf{MO\wedge_{R}} W)$ where $m=n-q$.
\end{theorem}

\begin{proof} The Pontrjagin-Thom construction establishes an
isomorphism of the bordism group $\M$ of manifolds with
$(\RO_W,\pi)$-structure in the stable tangent bundle and the homotopy group
$\pi_m(\mathbf{MO\wedge_R} W)$. Hence to prove the theorem it
suffices to show that the homomorphism $\psi$ is an isomorphism.

In what follows, in the notation of pullbacks of vector bundles we will occasionally suppress the symbols of maps that induce the bundles. For example, 
\[
J^k(EO_m, TW) \longrightarrow BO_m\times W
\]
will denote the formal $k$-jet bundle associated with the pullbacks of $EO_m$ and $TW$ with respect to the projections of $BO_m\times W$ onto the first and second factors respectively. 

\begin{lemma}\label{l:surj} The homomorphism $\psi$ is surjective.
\end{lemma}
\begin{proof} Given an element of $\M$, let an $(\RO_W, \pi)$-manifold $V_0$ be its representative. Suppose that the $(\RO_W, \pi)$-structure of $V_0$ is
represented by a continuous map $\alpha\co  V_0 \to [\RO_W]_{m+r}$ and a
fiberwise isomorphism of vector bundles
\[
  \tilde\alpha\co   TV_0\oplus r\varepsilon \longrightarrow \pi^*EO_{m+r}
\]
that covers $\alpha$ producing a commutative diagram
\[
\begin{CD}
  TV_0\oplus r\varepsilon @>\tilde \alpha>> \pi^*EO_{m+r}\\
  @VVV @VVV\\
  V_0  @>\alpha>> [\RO_{W}]_{m+r}.
\end{CD}
\]

We claim that the $(R_W, \pi)$-structure in the
stable tangent bundle of $V_0$ determines an $(R_{W\times L},
\pi_{W\times L})$-structure in the stable tangent bundle of
$V_0\times L$, where $L$ is a closed parallelizable manifold of dimension $r$ with a fixed trivialization of the tangent bundle. Indeed, the bundle
\[
[J^k_{W\times L}]_{m+r}\co = J^k(EO_{m+r}, T(W\times L))
\longrightarrow BO_{m+r}\times (W\times L)
\]
is canonically isomorphic to the bundle
\[
  [J^k_W]_{m+r}\times L\co = J^k(EO_{m+r}, TW\oplus r\varepsilon)\times
  L \longrightarrow (BO_{m+r}\times W)\times L
\]
under an isomorphism that takes the fiber
\[
J^k(EO_{m+r}|b, T(W\times L)|(w,l))  \qquad \textrm{over}\ (b,
(w,l))\in BO_{m+r}\times (W\times L)
\]
onto the fiber
\[
J^k(EO_{m+r}|b, (TW\oplus r\varepsilon)|w)\times\{l\}  \qquad
\textrm{over}\ ((b,w), l)\in (BO_{m+r}\times W)\times L.
\]
Consequently, since the basis of the differential relation for
$[R_{W\times L}]_{m+r}$ coincides with the basis of the differential
relation for $[R_W]_{m+r}$, the bundle
\[
[\RO_{W\times L}]_{m+r}\co = \RO(EO_{m+r}, T(W\times L))
\longrightarrow BO_{m+r}\times (W\times L)
\]
and the bundle
\[
  [\RO_W]_{m+r}\times L\co = \RO(EO_{m+r}, TW\oplus r\varepsilon)\times
  L \longrightarrow (BO_{m+r}\times W)\times L
\]
are canonically isomorphic. Furthermore, there is a canonical isomorphism of vector bundles
\[
\begin{CD}
  \pi^*EO_{m+r}\times L @>>> \pi_{W\times L}^*EO_{m+r}\\
  @VVV @VVV\\
  [\RO_W]_{m+r}\times L @>>> [\RO_{W\times L}]_{m+r},
\end{CD}
\]
which we precompose with the composition
\[
   T(V_0\times L) \longrightarrow (TV_0\oplus r\varepsilon)\times L \longrightarrow \pi^*EO_{m+r}\times L
\]
of a canonical isomorphism and the bundle map $\tilde\alpha\times\id_{L}$, where $\id_{L}$ is the identity map of $L$,  in order to obtain a desired $(\RO_{W\times L}, \pi_{W\times L})$-structure
\[
\begin{CD}
  T(V_0\times L) @>>> \pi^*_{W\times L}EO_{m+r}\\
  @VVV @VVV\\
  V_0\times L  @>>> [\RO_{W\times L}]_{m+r}
\end{CD}
\]
on the stable tangent bundle of $V_0\times L$. 

Similarly, the triple $(V_0, \alpha, \tilde\alpha)$ determines a fiberwise isomorphism of bundles 
\begin{equation}\label{eq:j1}
   J^k(T(V_0\times L), T(W\times L)) \longrightarrow [J^k_{W\times
   L}]_{m+r}.
\end{equation}
Let $\hat X\to V_0\times L$ be the trivial bundle with fiber
$W\times L$. Again, by choosing Riemannian metrics on $V_0\times L$ and $W\times L$, we
fix an isomorphism
\begin{equation}\label{eq:j2}
  J^k(\hat X)\longrightarrow J^k(T(V_0\times L), T(W\times L))
\end{equation}
and denote its composition with the bundle map (\ref{eq:j1}) by
\[
\jet\co  J^k(\hat X) \longrightarrow [J^k_{W\times
   L}]_{m+r}.
\]
By the h-principle, there is a homotopy $H=H_{\tau}$, with $\tau\in
[0,1]$, of the map
\[
H_0=\alpha\times\id_{L}\co  V_0\times L \longrightarrow
[R_W]_{m+r}\times L = [R_{W\times L}]_{m+r}
\]
to a section $H_1$ such that
$H_1=\jet\circ j^kh$ for some genuine smooth $\mathbf{R}$-map
\[
   h\co  V_0\times L \longrightarrow W\times L.
\]
By the Sard
Lemma applied to the composition of $h$ and the projection of
$W\times L$ onto the second factor, there is a value $pt\in L$
for which the map $h$ is transversal to the copy $W\times
\{pt\}$ of $W$. In particular the subset
$V_1\co =h^{-1}(W\times\{pt\})$ is a compact submanifold of $V_0\times L$.
Let $h_1\co  V_1\to W$ denote the restriction $h|V_1$ composed with
the identification $W\times \{pt\}\to W$. We will show that $h_1$
satisfies the differential relation $\mathbf{R}$, and $\psi([h_1])$
belongs to the class represented by $(V_0, \alpha, \tilde\alpha)$.

\begin{lemma}\label{p:1} The map $h_1\co  V_1\to W$ is a solution to
$\mathbf{R}$.
\end{lemma}
\begin{proof} For a point $x\in V_1$, by the Inverse
Function Theorem, there are neighborhoods $U(x)$ of $x$ in
$V_0\times L$ and $U(h(x))$ of $h(x)$ in $W\times L$ with
coordinates $(s_1,...,s_r,v_1,...,v_m)$ in $U(x)$ and $(\tilde
s_1,...,\tilde s_r, w_1,...,w_n)$ in $U(h(x))$ such that
$(0,\dots,0,v_1,...,v_m)$ are coordinates in $V_1\cap U(x)$, $(0,\dots,0,w_1,...,w_n)$
are coordinates in $(W\times\{pt\})\cap U(h(x))$, and the mapping
$h|U(x)$ has the form
\[
   \tilde s_i=s_i, \qquad i=1,...,r,
\]
\[
   w_j=w_j(s_1,...,s_r,v_1,...,v_m), \qquad j=1,...,n.
\]
We note that in the chosen coordinates the mapping $h|V_1\cap
U(x)$ has the form
\[
   w_j=w_j(0,...,0,v_1,...,v_m), \qquad j=1,...,n.
\]
Hence the local ring of the germ $h$ at $x$ is isomorphic to the
local ring of the germ $h|V_1$ at $x$. If $h|V_1$ is not a solution to a $\mathbf{K}$-invariant relation $\mathbf{R}$ in a neighborhood of $x$, then, since 
\[
   ext^{m}_{m+r}(J^k(\R^m, \R^{n})\setminus R_{m}) \subset J^k(\R^{m+r}, \R^{n+r})\setminus R_{m+r},
\] 
we conclude that $h$ is not an $\mathbf{R}$-map, which is a contradiction. Thus $h|V_1$ is a solution to $\mathbf{R}$.
 \end{proof}

It remains to show that $\psi([h_1])$ is the class
represented by $(V_0, \alpha, \tilde\alpha)$. 

To simplify the notation we will occasionally use the same symbol both for a map and its restrictions. In Lemma~\ref{l:10.4} we will chase the diagram
\[
\begin{CD}
[R_W]_m @>\subset>> [R_W]_{m+r} @>i_L>> [R_W]_{m+r}\times L @>\cong>>[R_{W\times L}]_{m+r} \\
@VVV @VVV @VVV @VVV \\
BO_m @>\subset>> BO_{m+r}@>>> BO_{m+r}@>>> BO_{m+r},
\end{CD}
\]
where the horizontal map $i_L$ is the inclusion of the $pt$-th slice $[R_{W}]_{m+r}\times \{pt\}$. 

\begin{lemma} \label{l:10.4}
The bordism class $\psi([h_1])$ is represented by an $(R_W, \pi)$-manifold $(V_1, \beta, \tilde\beta)$, where $\beta=H_1|V_1$ and $\tilde\beta$ is a bundle map covering $\beta$. 
\end{lemma}
\begin{proof} Suppose that the construction in the definition of $\psi$ applied to $h_1$ yields an $(R_W, \pi)$-manifold $(V_1, \gamma_h, \tilde \gamma_h)$ representing $\psi([h_1])$, with 
\[
     \gamma_h\co  V_1\longrightarrow [R_W]_m.
\] 
We recall that the triple $(V_1, \gamma_h, \tilde\gamma_h)$ is completely determined by the choices of Riemannian metrics on $V_1$ and $W$ and a bundle map 
\[
   \tau\co  TV_1\longrightarrow EO_m. 
\]
We pass to the $r$-th suspension of $(V_1, \gamma_h, \tilde\gamma_h)$ and denote it by the same symbols so that now for example 
\[
   \gamma_h\co  V_1\longrightarrow [R_W]_{m+r}. 
\]
On the other hand, the map	
\[
    \beta=H_1|V_1\colon V_1\longrightarrow [R_{W\times L}]_{m+r}
\]
factors through a map 
\[
    \gamma_H\colon  V_1\longrightarrow [R_W]_{m+r},
\]
while the map $\tilde\beta$ factors through a map $\tilde\gamma_H$ covering $\gamma_H$. We need to show that the $(R_W, \pi)$-manifolds $(V_1, \gamma_h, \tilde\gamma_h)$ and $(V_1, \gamma_H, \tilde\gamma_H)$ represent the same bordism class. In fact we will construct a homotopy of $(\gamma_H, \tilde\gamma_H)$ to $(\gamma_h, \tilde\gamma_h)$.  

In view of the Riemannian metric on $L$ and a fixed trivialization of $TL$, we may identify a neighborhood $U$ of $0$ in the vector space $\R^r$ with a neighborhood of $pt$ in $L$. Furthermore, we may identify a tubular neighborhood of $V_1$ in $V_0\times L$ with $V_1\times U$ so that the restriction of $h$ to $V_1\times U$ is of the form 
\[
       V_1\times U\longrightarrow W\times U,
\]
\[
       (v,u)\mapsto (\bar{h}_1(v,u), u)
\]
for some map $\bar{h}_1$. This map is fiberwise (over $U$) homotopic to the map $h_1\times \id_U$ through $\mathbf{R}$-maps. Consequently, we may assume that $h|V_1\times U$ coincides with $h_1\times \id_U$. 

Let $\hat Y\to V_1\times U$ denote the trivial bundle with fiber $W\times U$. Then there is an obvious inclusion $\hat Y\subset \hat X$ which in its turn determines an inclusion $J^k(\hat Y)\subset J^k(\hat X)$. Let us recall that the restriction $H_1|V_1\times U$ is given by the composition of 
\[
j^kh\colon V_1\times U\longrightarrow J^k(\hat{Y})\subset J^k(\hat{X}) 
\] 
and 
\[
   \jet\colon J^k(\hat Y)\stackrel{(\ref{eq:j2})}\longrightarrow J^k(T(V_1\times U), T(W\times U))\stackrel{(\ref{eq:j1})}\longrightarrow [J^k_{W\times U}]_{m+r}.
\]
To begin with we modify the Riemannian metrics on $V_1\times U$ and $W\times U$ by homotopy to Riemannian metrics given by products of those already fixed on $V_1$, $W$ and $U$. This determines a homotopy of the map (\ref{eq:j2}), and therefore a homotopy of $H_1|V_1\times U$ and $(\gamma_H, \tilde\gamma_H)$. 

Next we recall that the map (\ref{eq:j1}) is determined by a bundle map 
\[ T(V_1\times U)=TV_1\times TU \longrightarrow EO_{m+r}.
\]
We modify it by homotopy to a bundle map given by the composition
\[
\begin{CD}
   TV_1\times TU @>>> TV_1\oplus \varepsilon^r @>>> EO_m\oplus \varepsilon^r @>>> EO_{m+r}\\
   @VVV @VVV @VVV @VVV \\
   V_1\times U @>>> V_1 @>>> BO_m @>>> BO_{m+r},
\end{CD}
\]
where the left horizontal maps are obvious projections along the factor $U$, the middle horizontal maps are determined by $\tau$, and the right horizontal maps are inclusions. 
This determines a homotopy of (\ref{eq:j1}) and a further homotopy of $H_1|V_1\times U$ and $(\gamma_H, \tilde\gamma_H)$. 

Let us observe now that the obtained pair $(\gamma_H, \tilde\gamma_H)$ coincides with $(\gamma_h, \tilde\gamma_h)$. 
\end{proof}

We may perturb the homotopy $H$ relative to $H_0$ and $H_1$ so that $H$
becomes transversal to the submanifold
\begin{equation}\label{eq:4}
[\RO_{W}]_{m+r}\times \{pt\}\times [0,1] \quad \textrm{in} \quad
[\RO_W]_{m+r}\times L \times [0,1].
\end{equation}
Then the inverse image $V$ of
$[\RO_W]_{m+r}\times\{pt\}\times [0,1]$ under the map $H$ is a compact
submanifold of $V_0\times L\times [0,1]$ with boundary $\partial V$ that consists of two parts
\[
   V_0\times \{pt\}\times \{0\}\quad \textrm{and} \quad V_1 \times
   \{1\}.
\]

Now it suffices to prove the assertion that the composition

\begin{eqnarray}
   V&\stackrel{H|V}\longrightarrow& [\RO_W]_{m+r}\times\{pt\}\times [0,1] \nonumber\\
   &\stackrel{\subset}\longrightarrow& [\RO_W]_{m+r+1}\times\{pt\}\times [0,1] \nonumber\\
   &\longrightarrow& [\RO_W]_{m+r+1} \nonumber
\end{eqnarray}
can be covered by a
bundle map 
\[
 TV\oplus r\varepsilon \longrightarrow \pi^*EO_{m+r+1}
\]
that leads to an $(\RO_W,\pi)$-bordism between the $1$-step stabilizations of $(V_0, \alpha, \tilde\alpha)$ and $(V_1, \beta, \tilde \beta)$.

To verify the assertion, let us observe that the $1$-step stabilization of the bundle map $\tilde\alpha$ extends to a bundle map $\tilde H$, 
\[
\begin{CD}
     T(V_0\times L\times [0,1]) @>\tilde H >> (\pi_{W\times L} \circ\pi_{1,2})^*EO_{m+r+1} \\
@VVV @VVV \\
V_0\times L\times [0,1] @>>> [R_W]_{m+r+1}\times L\times [0,1],
\end{CD}
\]
covering the composition
\[
   H'\colon V_0\times L\times [0,1] \stackrel{H}\longrightarrow [R_W]_{m+r}\times L\times [0,1] \stackrel{\subset}\longrightarrow [R_W]_{m+r+1}\times L\times [0,1],
\]
where $\pi_{1,2}$ is the projection of $[R_{W}]_{m+r+1}\times L\times [0,1]$ onto the product of the first two factors. 

The normal bundle $\nu$ of $V$ in $V_0\times L\times [0,1]$ is trivial as it is isomorphic to the pullback via $H$ of the trivial normal bundle of the submanifold (\ref{eq:4}). Consequently the restriction of the bundle map $\tilde H$ to $V$,
\[
\begin{CD}
     TV\oplus\nu @>\tilde{H}|V >> (\pi \circ\pi_{1,2})^*EO_{m+r+1} @>=>> \pi_1^*EO_{m+r+1} \\
@VVV @VVV @VVV\\
V @>H'>> [R_W]_{m+r+1}\times \{pt\}\times [0,1] @>=>> [R_W]_{m+r+1}\times[0,1],
\end{CD}
\]
where $\pi_1$ is the projection of $[R_W]_{m+r}\times [0,1]$ onto the first factor,
leads to an $(R_W, \pi)$-bordism of the $1$-step stabilizations of the $(R_W, \pi)$-manifolds $(V_0, \alpha, \tilde\alpha)$ and $(V_1, \beta, \tilde\beta)$. 

 \end{proof}

\begin{lemma} The homomorphism $\psi$ is injective.
\end{lemma}
\begin{proof} Let $f\co V_0\to W_0$ be a map representing an
element $[f]\in SB(W_0)$ in the kernel of the homomorphism
$\psi$. Then $f$ determines an $(\RO_{W_0},\pi)$-structure in the
stable tangent bundle of $V_0$ such that for a sufficiently big
$r$, the structure map $\tilde \alpha_0$ in the commutative
diagram
\[
\begin{CD}
   TV_0\oplus r\varepsilon @>\tilde\alpha_0>>
   \pi^*EO_{m+r}\\
   @VVV @VVV \\
   V_0 @>\alpha_0>> [\RO_{W_0}]_{m+r}
\end{CD}
\]
extends to a structure map $\tilde\alpha$ in the commutative
diagram
\[
\begin{CD}
   TV\oplus (r-1)\varepsilon @>\tilde\alpha>>
   \pi^*EO_{m+r}\\
   @VVV @VVV \\
   V @>\alpha>> [\RO_{W}]_{m+r}
\end{CD}
\]
where $W=W_0\times [0,1]$, and $V$ is a compact $(m+1)$-dimensional
manifold with boundary $\partial V=V_0$.

Let $\hat X\to V\times L'$ be the trivial bundle with fiber
$W\times L'$, where $L'$ is a parallelizable closed manifold of dimension $r-1$ with a fixed trivialization of the tangent bundle. As in the proof of Lemma~\ref{l:surj}, we use
the map $\tilde\alpha$ to construct a commutative diagram
\[
\begin{CD}
   T(V\times L') @>\tilde\delta>> \pi_{W\times L'}^*EO_{m+r}\\
   @VVV @VVV \\
   V\times L' @>\delta>> [\RO_{W\times L'}]_{m+r}
\end{CD}
\]
that defines an $(\RO_{W\times L'}, \pi_{W\times L'})$-structure in the stable tangent bundle of $V\times L'$. Also,
as in the proof of Lemma~\ref{l:surj}, we define a mapping
\[
   \jet\co  J^k(\hat X)\longrightarrow [J^k_{W\times L'}]_{m+r}.
\]
Let $C\approx V_0\times [0,a)$ for some $0<a\ll 1$ be a small collar neighborhood of
$V_0$ in $V$. We may assume that in $C\times L'$ the
structure map $\delta$ coincides with $\jet\circ
j^k[f\times\id_{[0,a)}\times \id_{L'}]$. Then, by the relative
h-principle applied to the pair $(V\times L', C\times
L')$, there is a homotopy of $\delta$ to a map $V\times
L'\to [\RO_{W\times L'}]_{m+r}$ given by the composition
$\jet\circ j^k h$ for some genuine solution $h\co V\times L'\to W\times
L'$ of $\mathbf{R}$.

Again, as in the proof of Lemma~\ref{l:surj}, there is a regular
value $pt$ of the composition of $h$ and the projection of
$W\times L'$ onto the second factor. For a compact submanifold $V_1$ defined as
$h^{-1}(W\times \{pt\})$ we can show that
\[
   h|V_1\co  V_1\longrightarrow W\times \{pt\},
\]
is a solution of $\mathbf{R}$. Since
$h|\partial V_1$ can be identified with $f$, this implies that the map represents the trivial
element in $SB(W)$.

 \end{proof}

The proof of Theorem~\ref{t:iso} is complete.
 \end{proof}

%

\section{Proof of the weak bordism principle}\label{s:11}

Let $\mathbf{Diff}_n$ denote the category of smooth (possibly non-compact) $n$-mani\-folds without boundary and smooth (possibly non-proper) embeddings. 
In this section we consider the covariant functor 
\[
          \EB(\cdot)\co  \mathbf{Diff}_n \longrightarrow \mathbf{AG} 
\]
(see section~\ref{s:5}) that takes a manifold $W$ of dimension $n$ onto the bordism group $\EB(W)$ of $\mathbf{R}$-maps of $(n-q)$-dimensional manifolds into $W$, and  
prove the weak bordism principle (see Definition~\ref{def:2}) for $\EB(\cdot)$
with respect to the category $\Top^2\downarrow BO_n$ of pairs of topological spaces over the classifying space $BO_n$ of vector $n$-bundles. 

To begin with, let us observe that the definition of the bordism group $\mathbf{MO_R}(W)$ of stable formal $\mathbf{R}$-maps into an $n$-manifold $W$ can be easily extended to that of the bordism group of stable formal $\mathbf{R}$-maps into a pair $(X, A; \varphi)$ of topological spaces over $BO_n$. 

Indeed, a $\mathbf{K}$-invariant sequence $\mathbf{R}=\{R_s(q)\}$ of differential relation basis $R_s\subset J^k(\R^s, \R^{q+s})$, with $s\ge 0$ and $s\ge n-q$,  leads to a sequence of formal differential relations $[R_X]_s$ in the formal $k$-jet bundles $[J_X]_s$ associated to the pair of bundles
\[
     (\pi_1^*EO_s, (\varphi\circ \pi_2)^*EO_n\oplus t\varepsilon)   \qquad \textrm{over} \qquad BO_s\times X,
\]
where $\pi_1$ and $\pi_2$ denote the projections of $BO_s\times X$ onto the first and second factors respectively, and $t=q+s-n$. Following Definition~\ref{def:3}, we define the bordism group $\mathbf{MO_{R}}(X; \varphi)$ of stable formal $\mathbf{R}$-maps into $(X; \varphi)$ as the $m$-th bordism group of maps of manifolds with $(R_X, \pi)$-structure in the stable tangent bundle, where $m=n-q$, and $(R_X, \pi)$ is the $(B,f)$-sequence given by commutative diagrams as in Definition~\ref{def:3}. In other words, $\mathbf{MO_R}(X; \varphi)$ is the $m$-th homotopy group of the Thom spectrum $\mathbf{MO\wedge_R}(X; \varphi)$ associated to the tangent $(R_X, \pi)$-structures. 
Finally we define the group
$\mathbf{MO_{R}}(X, A; \varphi)$ as the relative $m$-th homotopy group of the pair of $\mathbf{MO\wedge_R}(X; \varphi)$ and its subspectrum $\mathbf{MO\wedge_R}(A; \varphi|A)$. 

More generally we define the group 
$\mathbf{MO}_{\mathbf{R}\, l}(X, A; \varphi)$ as the relative $l$-th homotopy group of the pair of $\mathbf{MO\wedge_R}(X; \varphi)$ and $\mathbf{MO\wedge_R}(A; \varphi|A)$. Furthermore, for each $l\ge 0$, we define a functor $\mathbf{MO}_{\mathbf{R}\, l}$ on the category $\Top^2 \downarrow BO_n$ so that 
\[
             \mathbf{MO}_{\mathbf{R}\, l}\colon \Top^2 \downarrow BO_n \longrightarrow \mathbf{AG}
\]
is given by 
\[
        \mathbf{MO}_{\mathbf{R}\, l}\colon (X, A, \varphi)\  \mapsto\  \mathbf{MO}_{\mathbf{R}\, l}(X, A, \varphi)
\]
and for each morphism 
\[
    i\co  (X_1, A_1, \varphi_1) \longrightarrow (X_2, A_2, \varphi_2)
\]
of pairs of topological spaces over $BO_n$, the homomorphism
\[
       \mathbf{MO}_{\mathbf{R}\, l}(i)\colon \mathbf{MO}_{\mathbf{R}\, l}(X_1, A_1, \varphi_1) \longrightarrow \mathbf{MO}_{\mathbf{R}\, l}(X_2, A_2, \varphi_2)
\]
is induced by the obvious map of pairs of spectra
\[
     (\mathbf{MO\wedge_R}(X_1; \varphi_1), \mathbf{MO\wedge_R}(A_1; \varphi_1)) \longrightarrow      (\mathbf{MO\wedge_R}(X_2; \varphi_2), \mathbf{MO\wedge_R}(A_2; \varphi_2)).
\]

In view of the Blakers-Massey theorem, it is easily verified that $\mathbf{MO}^{*}_\mathbf{R}(\cdot)$ determines a homology theory on the category of pairs of topological spaces over $BO_n$. 

To complete the proof of the weak bordism principle for the functor $\EB(\cdot)$, it suffices to observe that Theorem~\ref{t:iso} implies that $\EB(\cdot)$ is naturally equivalent to the functor $\mathbf{MO_R}(\cdot)$ precomposed with the functor 
\[
     G\co  \mathbf{Diff}_n \longrightarrow \Top^2 \downarrow BO_n
\]
such that  
\begin{itemize}
\item $G$ takes a manifold $W$ of dimension $n$ onto the pair $(X; \varphi)$ of the topological space $X$ underlying $W$ and a map $\varphi\co  X\to BO_n$ classifying the vector bundle inherited from $TW$; and 
\item $G$ takes a smooth embedding 
\[
    i\co  W_1 \longrightarrow W_2
\]
in $\mathbf{Diff}_n$ onto the morphism 
\[
    G(i)\co  (X_1, \varphi_1)\longrightarrow (X_2, \varphi_2),
\] 
where $(X_i, \varphi_i)\co =G(W_i)$ for $i=1,2$,  that as a map of topological spaces is given by the composition 
\[
      X_1 \stackrel{\cong}\longrightarrow W_1 \stackrel{i}\longrightarrow W_2 \stackrel{\cong}\longrightarrow X_2.
\] 
\end{itemize}

\begin{theorem}\label{th:1/2} Under the assumptions of the Main Theorem, the covariant functor $\EB$ satisfies the weak bordism principle. 
\end{theorem}

\section{A contravariant functor dual to $\mathbf{MO_R}(\cdot)$}\label{s:13}

As it has been shown in the previous section, the functor $\mathbf{MO_R}(\cdot)$ associated to the functor 
\[
          \EB(\cdot)\co  \mathbf{Diff}_n \longrightarrow \mathbf{AG} 
\]
of bordism groups of maps satisfying a $\mathbf{K}$-invariant relation $\mathbf{R}$ is defined only on the category $\Top^2\downarrow BO_n$ of pairs of topological spaces over $BO_n$. In particular a number of effective methods of homotopy theory are not available for $\mathbf{MO_R}(\cdot)$. We will show, however, that there is a contravariant functor $H_{\mathbf{R}}^0(\cdot)$ dual to $\mathbf{MO_R}(\cdot)$ such that $H_{\mathbf{R}}^0(\cdot)$ is a cohomology functor in the classical sense, and 
there is a duality isomorphism
\[
    H^0_{\mathbf{R}}(W)\cong \mathbf{MO_R}(W)
\] 
for any closed manifold $W$ of dimension $n$. 

To define $H_{\mathbf{R}}^*(\cdot)$, let us recall that $EO_t\to BO_t$, with $t\ge 0$, denotes the universal vector $t$-bundle. There is a fiber bundle $\pi\co  S_t\to BO_t$ with the total space given by the space of $k$-jets of $\mathbf{R}$-maps
\begin{equation}\label{eq:5} 
   (\R^{t-q}, 0) \longrightarrow (EO_t|b, 0),  \qquad \textrm{with} \qquad b\in BO_t,
\end{equation}
and with projection $\pi$ taking the $k$-jet of a germ~(\ref{eq:5}) onto the point $b\in BO_t$. A construction similar to that in Definition~\ref{def:3} provides us with a $(B,f)$-sequence $(S, \pi)$ given by the commutative diagrams 
\[
\begin{CD}  S_t @>>> S_{t+1} \\
    @V\pi VV @V\pi VV \\
    BO_t @>>> BO_{t+1}
\end{CD}
\]
where the horizontal maps are the canonical inclusions as in Definition~\ref{def:3} (cf. the definition of $\mathbf{B}_{\J}$ in section~\ref{s:1.1} where $k=\infty$).  

\begin{definition} The Thom spectrum associated to the normal $(S, \pi)$-struc\-tures defines a generalized cohomology theory (in the classical sense). Its contravariant functors are denoted by $\{H^*_{\mathbf{R}}(\cdot)\}$. 
\end{definition}

For a manifold $W$ of dimension $n'$, the group $H_{\mathbf{R}}^r(W)$ is the $r$-th cobordism group of maps into $W$ of manifolds with an $(S, \pi)$-structure in the stable normal bundle. In particular, every element of $H^r_{\mathbf{R}}(W)$ is represented by a $4$-tuple $(V, i, \alpha, \tilde\alpha)$ of 

\begin{itemize}
\item a manifold $V$ of dimension $m=n'-r-q$, 
\item an embedding $i\co  V\to \R^{t-r-q}\times W$ for sufficiently big $t$ such that the composition of $i$ and the projection $\R^{t-r-q}\times W\to W$ onto the second factor is proper, 
\item a continuous map $\alpha\co  V\to S_t$, and
\item a fiberwise isomorphism $\tilde\alpha\co  \nu\to \pi^*EO_t$ covering $\alpha$, where $\nu$ is the normal vector bundle of $V$ induced by the embedding $i$. 
\end{itemize}

\begin{theorem}\label{th:2} Under the assumptions of the Main Theorem there is a canonical isomorphism  $H^0_{\mathbf{R}}(W)\to \EC(W)$ for any closed manifold $W$ of dimension $n$. 
\end{theorem}
\begin{proof} Let $(V, i, \alpha, \tilde\alpha)$ be a $4$-tuple representing a given element of $H^0_{\mathbf{R}}(W)$. We note that since $W$ is closed, the manifold $V$ is closed as well.  Let
\begin{equation}\label{eq:7}
      J^k((t-q)\varepsilon, \nu) \longrightarrow V
\end{equation}
be the formal $k$-jet bundle over $V$ associated to the pair of bundles over $V$, namely, $(t-q)\varepsilon$ and the normal vector $t$-bundle $\nu$ induced by $i$. We define a section $s$ of the bundle~(\ref{eq:7})
by taking a point $v\in V$ onto the $k$-jet represented by the composition of a map germ~(\ref{eq:5}) representing $\alpha(v)$ and the composition of two isomorphisms
\[
    EO_t|b \xrightarrow{(\pi|\alpha(v))^*} ((\pi|\alpha(v))^*)(EO_t|b) \xrightarrow{(\tilde\alpha|v)^{-1}} \nu|v,
\]
where $b=\pi(\alpha(v))$.
There is a canonical inclusion 
\begin{equation}\label{eq:6}
   J^k((t-q)\varepsilon, \nu) \longrightarrow J^k(TV\oplus  (t-q)\varepsilon, TV \oplus \nu)
\end{equation}
of total spaces of $k$-jet bundles over $V$ that takes the $k$-jet at $v\in V$ represented by a germ $f$ onto the $k$-jet at $v$ represented by the germ $\id_{TV|v} \times f$. Let 
\[
    \tilde{s}\co  V\longrightarrow J^k(TV\oplus (t-q)\varepsilon, TV\oplus \nu)
\]
 be the composition of $s$ and the canonical inclusion~(\ref{eq:6}). Let us observe that $TV\oplus \nu$ is canonically isomorphic to the pullback of $(t-q)\varepsilon \oplus TW$ with respect to the composition 
\[
    V\longrightarrow  \R^{t-q}\times W \longrightarrow W
\]
of the embedding $i$ and the projection of $\R^{t-q}\times W$ onto the second factor. Consequently, $\tilde s$ determines an $(R_W, \pi)$-structure on the stable tangent bundle of $V$, which in its turn, by Theorem~\ref{t:iso},
determines an element of $\EC(W)\cong \EB(W)$. It is easily verified that the obtained element in $\EC(W)$ does not depend on the choice of the representative of $H^0_{\mathbf{R}}(W)$. Thus, there is a canonical map
\[
     \xi_W\co  H^0_{\mathbf{R}}(W) \longrightarrow \EC(W)
\]
which is, in fact, a homomorphism. A similar construction provides us with a
homomorphism 
\[
   \xi^{-1}_W\co  \EC(W) \longrightarrow H^0_{\mathbf{R}}(W)
\]  
inverse to $\xi_W$. Thus the homomorphism $\xi_W$ is an isomorphism. 
 \end{proof}

It is easily verified that the homomorphisms $\xi^{-1}_W$ in Theorem~\ref{th:2} determine a natural transformation of functors on the category of closed manifolds of dimension $n$. 

Now the Main Theorem readily follows from Theorem~\ref{th:2}.

\begin{theorem} Under the assumptions of the Main Theorem, the contravariant functor $\EC$ restricted to the full subcategory of closed manifolds of dimension $n$ satisfies the b-principle. 
\end{theorem}

\section{Applications}\label{s:app}

As has been mentioned, for a sequence $\mathbf{R}$ of $\K$-invariant bases, the set of right bordism classes of $\mathbf{R}$-maps into a manifold $W$ forms a semigroup, which is a priori not a group. On the other hand, since $\mathbf{MO_R}(W)$ is always a group, Theorem~\ref{t:iso} implies that this semigroup is often a group.  

\begin{corollary}\label{c:1} If a sequence $\mathbf{R}$ of bases satisfies the conditions of Theorem~\ref{t:iso}, then the right bordism classes of $\mathbf{R}$-maps into a manifold $W$ form a group.  
\end{corollary}

We note that the assumptions in Corollary~\ref{c:1} are necessary as, for example, the right bordism classes of submersions into a manifold $W$ do not form a group. Corollary~\ref{c:1} was also independently observed by Sz\H{u}cs~\cite{Sz1} in the negative dimension case and by Ando~\cite{An4} in the general case. 

%

Furthermore, the weak b-principle proved for the functor $\EB$ allows us to use the machinery of cohomology theory, notably spectral sequences and theorems on homology of infinite loop spaces, in order to carry out explicit computations ~\cite{Sad4} (see also~\cite{Sz1}). 

The class of differential relations for which our theorem applies is large. For example, as we will see next, it contains almost all Thom-Boardman differential relations. 

\subsection{Maps with prescribed Thom-Boardman singularities} 

Thom-Boardman singularities $\Sigma_{\I}$ are singularities of smooth maps, each indexed by a sequence $\I=(i_1,...,i_l)$ of integers. For a
definition and properties of Thom-Boardman singularities we refer the reader to the original paper of Boardman \cite{Bo}.  

\begin{definition} Given a sequence $\I=(i_1,\dots, i_l)$ of integers, the {\it Thom-Boardman basis} $R^{\I}_s\subset J^k(\R^{s}, \R^{q+s})$ is defined to be the complement to the set $\cup\Sigma_{\I'}$ where the union is taken over $\I'>\I$ with respect to the lexicographic order. It is well-known that the relations $R^{\I}_s$, $s\ge 0$, form a sequence, which we will denote by $\mathbf{R}^{\I}$, of open $\K$-invariant bases. 
\end{definition}

It follows that solutions of
differential relations associated with the basis $R_s^{\I}$ are the maps without Thom-Boardman singularities $\Sigma_{\I'}$ with index $\I'>\I$.

The homotopy principle for Thom-Boardman bases $R_s^{\I}\subset J^k(\R^s, \R^{q+s})$, $q\le 0$, was proved by A.~du Plessis~\cite{Pl1} and a version of the h-principle in the case of the symbol $(1-q,0)$ was proved by Eliashberg~\cite{El1},\cite{El2} (see also Ando \cite{An10}), and then, in the case of an arbitrary symbol $\I$ with $\I>(1-q,0)$, using the Eliashberg h-principle, by Ando~\cite{An7}. Consequently, we derive the  bordism principle for Thom-Boardman differential relations. 

\begin{corollary}\label{c:TB} Suppose that either $q>0$ or $q\le 0$ and $\I\ge (1-q, 0)$. Then for the Thom-Boardman sequence of bases $R_s^{\I}\subset J^k(\R^s, \R^{q+s})$ the weak bordism principle holds for the corresponding functor $\EC$ on the category of closed manifolds and the corresponding functor $\EB$ on the category of manifolds. 
\end{corollary}


\subsection{Maps with additional structure}
It is easy to extend the results to the case of bordism groups of manifolds with an additional structure, e.g., oriented manifolds, manifolds with almost complex structure or spin manifolds. We refer the reader to the paper \cite{El} for the necessary adjustments in the case of Lagrangian and Legendrian immersions.  
 
\subsection{Kazarian conjecture} Under the assumptions of the Main Theorem, for oriented $\J$-cobordism groups we can deduce the Kazarian conjecture~\cite{Sz1}, which relates rational cohomology groups $H^*(B_{\J}; \mathbb{Q})$ of the classifying space $B_{\J}$ of $\J$-cobordism groups   with the cohomology groups of $\mathbf{S}(\J)=\lim S_t(\J)$, provided of course that $B_{\J}$ exists. Namely it was conjectured that the classifying space of $\J$-cobordism groups is rationally equivalent to the space 
\[
\lim_{t\to\infty}\ \Omega^{t+d}S^{t+d} S_t(\J).
\]
On the other hand,  the latter space is rationally equivalent to $\Omega^{\infty}\B_{\J}$. Since $\Omega^{\infty}\B_{\J}$ is itself a classifying space of $\J$-cobordism groups (by the Main Theorem), we establish the Kazarian conjecture. In the negative dimension case the Kazarian conjecture is proved by Sz\H{u}cs \cite{Sz1}.

\subsection{Prospective applications} In view of the Rim\'anyi-Sz\H{u}cs construction $B_{\J}$ and the Kazarian construction $K_{\J}$ of classifying spaces, the Main Theorem can be used to study diffeomorphism groups of smooth manifolds (see Remark~\ref{r:1.10} and the discussion before it).

\section{Obstruction spectrum}\label{s:15}

Let $\mathbf{R}$ be a sequence of differential relation bases $\RO_s\subset  J^k(\R^{s},\R^{q+s})$, with $s\ge 0$, for which Theorem~\ref{th:2} holds. In
this section we define\footnote{For morphisms of vector bundles
similar obstructions have been defined by Koschorke in \cite{Ko}. For
negative dimension mappings without certain singularities similar obstructions are
defined and studied by Sz\H ucs in \cite{Sz}, \cite{Sz1}. For an equivariant
cohomology version of these obstructions we refer the reader to
\cite{FR}.} a complete obstruction to the existence of a cobordism of
a given smooth map to an $\mathbf{R}$-map.

Let $\mathbf{R'}=\{R'_s\}$ be the sequence of full differential relation bases $R'_s=J^k(\R^s, \R^{q+s})$, with $s\ge 0$. Then the inclusions $R_s\subset R'_s$ lead to an inclusion $i$ of the Thom spectrum of $H^*_{\mathbf{R}}$ into the Thom spectrum of the cobordism group $H^{*-q}_{\mathbf{R'}}=\mathcal{N}^{*}$. Let $\mathcal{O}^{*}_{\mathbf{R}}$ denote the cofiber of $i$, the {\it obstruction spectrum}. Then for a manifold $W$ there is a long
exact sequence of groups 
\[
  \cdots \longrightarrow H^{*-q}_{\mathbf{R}}(W) \longrightarrow \mathcal{N}^{*}(W)
   \stackrel{j}\longrightarrow \mathcal{O}_{\mathbf{R}}^{*}(W) \longrightarrow
  H^{*+1-q}_{\mathbf{R}}(W) \longrightarrow \cdots.
\]
A proper smooth map $f\co V\to W$ of dimension $-q$ represents an element $[f]$ in the cobordism group  $\mathcal{N}^q(W)$. We define the obstruction class $o(f,\mathbf{R})\in \mathcal{O}^q_{\mathbf{R}}(W)$ as the image of $[f]$ 
with respect to the homomorphism $j$.  

\begin{theorem}\label{th:obst} Let $f\co V \to W$ be a map of dimension $-q$ of closed manifolds. Then $o(f,\mathbf{R})\in \mathcal{O}^q_{\mathbf{R}}(W)$ is a complete obstruction to the existence of a cobordism of $f$ to an $\mathbf{R}$-map.
\end{theorem}

In the Habilitation thesis~\cite{Ka3} (see also~\cite{Sz}), Kazarian presented obstructions to the existence of a cobordism of a given map to a map without certain multi-singularities as cohomology operations with values in singular cohomology groups. Similarly, we may view the map $j\co  \mathcal{N}^*(W) \to \mathcal{O}_{\mathbf{R}}^*(W)$ as a cohomology operation with values in an extraordinary cohomology theory. Thus if $\mathbf{R}$ satisfies the assumptions of Theorem~\ref{th:obst}, then the complete obstruction to the existence of a bordism to an $\mathbf{R}$-map is given by a cohomology operation with values in an extraordinary cohomology theory.


\begin{thebibliography}{99}
\bibitem{An1} Y.~Ando, Folding maps and the surgery theory on
manifolds, J. Math. Soc. Japan, 53 (2001), 357--382.
\bibitem{An10} Y.~Ando, Existence theorems of fold-maps, Japan. J. Math., 30 (2004), 29--73.
\bibitem{An5} Y.~Ando, The homotopy principle for maps with singularities of given $\mathcal K$-invariant class,  J. Math. Soc. Japan,  59 (2007), 557--582.
\bibitem{An6} Y.~Ando, Smooth maps with singularities of bounded $\mathcal{K}$-codimensions, arXiv:math/0704.0115v2.
\bibitem{An7} Y.~Ando, A homotopy principle for maps with prescribed Thom-Boardman singularities, Trans. Amer. Math. Soc., 359 (2007), 489--515.  
\bibitem{An4} Y.~Ando, Cobordisms of maps with singularities of a given class, Alg. Geom. Topol., 8 (2008), 1989--2029. 
\bibitem{Bo}  J.~M.~Boardman, Singularities of differentiable
maps, Inst. Hautes \'Etud. Sci. Publ. Math.,  33 (1967), 21--57.
\bibitem{Ch} D.~S.~Chess, Singularity theory and configuration space models of $\Omega^{n}S^n$ of nonconnected spaces, Topology Appl., 25 (1987), 313--338.
\bibitem {El1} Y.~Eliashberg, On singularities of folding type,
Math. USSR, Izv., 4 (1970), 1119--1134.
\bibitem {El2} Y.~Eliashberg, Surgery of singularities of smooth mappings, Math. USSR, Izv., 6 (1972), 1302--1326.
\bibitem{El} Y.~Eliashberg, Cobordisme des solutions de
relations diff\'erentielles, South Rhone seminar on geometry. I,
Lyon, 1983, Trav. Cours, 1984, 17--31.
\bibitem{EM1} Y.~Eliashberg, N.~Mishachev, Wrinkling of smooth mappings III. Foliations of codimension greater than one,  Topol. Methods Nonlinear Anal., 11 (1998), 321--350.
\bibitem{EM} Y.~Eliashberg,  N.~Mishachev, Introduction to the
h-principle, Grad. Stud. Math., 48, AMS, Providence, RI,
2002.
\bibitem{FR} L.~M.~Feh\'er, R.~Rim\'anyi, Calculation of Thom
polynomials and other cohomological obstructions for group
actions, Real and Complex Singularities (S\~ao Carlos, 2002), Ed. T.Gaffney and M.Ruas, Contemp. Math., 354, AMS, Providence, RI (2004), 69--93. 
\bibitem{Fu} D.~B.~Fuks, Quillenization and bordism, Funkts. Anal. Prilozh., 8 (1974), 36--42; translation in Funct. Anal. Appl., 4 (1970), 143–-151.
\bibitem{EG} S.~Galatius, Y.~Eliashberg, Homotopy theory of compactified moduli space, Oberwolfach Report 13/2006, 761--767.
\bibitem{GMTW} S.~Galatius, I.~Madsen, U.~Tillmann, M.~Weiss, The homotopy type of the cobordism category, preprint, arXiv:math/0605249v2, 2007.
\bibitem{Gr71} M.~Gromov, A topological technique for the 
construction of solutions of differential equations and inequalities,
Proc. Int. Congress Math., 2 (1970), 221--225.
\bibitem {Gr} M.~Gromov, Partial differential relations,
Springer-Verlag, Berlin, Heidelberg, 1986.
\bibitem{EGr} M.~Gromov, Y.~Eliashberg, Removal of singularities of smooth mappings, Izv. Akad. Nauk SSSR, 35 (1971), 600--626; translation in Math. USSR Izvestija, 5 (1971), 615--639.
\bibitem{Hat} A.~Hatcher, Algebraic Topology, Cambridge Univ. Press, Cambridge, 2002.
\bibitem{Ik} K.~Ikegami, Cobordism group of Morse functions on manifolds, Hiroshima Math. J., 34 (2004), 211--230.
\bibitem {IS} K.~Ikegami, O.~Saeki, Cobordism group of Morse
functions on surfaces, J. Math. Soc. Japan, 55 (2003), 1081--1094.
\bibitem{Ja} K.~J\"anich, Symmetry properties of singularities
of $C^{\infty}$-functions, Math. Ann., 238 (1978), 147--156.
\bibitem{Kal0} B.~Kalm\'ar, Cobordism group of Morse functions on unoriented surfaces, Kyushu J. Math., 59 (2005), 351--363.
\bibitem{Kal2} B.~Kalm\'ar, Fold cobordisms and stable homotopy groups, preprint, arXiv:math/0704.3147.
\bibitem{Ka3} M.~E.~Kazarian, Kharakteristicheskie classy v teorii 
osobennostej, Habilitation thesis, 2003. 
\bibitem{Ka1} M.~Kazarian, Multisingularities, cobordisms, and enumerative geometry (Russian), Uspekhi Mat. Nauk, 58 (2003), 29--88; English translation in Russian Math. Surv., 58 (2003), 665--724.
\bibitem{Ka2} M.~Kazarian, Thom polynomials, singularity theory and its applications; Edited by S.~Izumiya, G.~Ishikava, H.~Tokunaga, I.~Shimada and T.~Sano, Adv. Stud. Pure Math., 43 (2007), 85--135.
\bibitem{Koch} S.~O.~Kochman, Bordism, stable homotopy and Adams spectral sequences, Fields Inst. monographs, AMS, Providence, Rhode Island, 1996. 
\bibitem{Ko} U.~Koschorke, Vector fields and other vector bundle
morphisms--a singularity approach, Lect. Notes Math., 847,
Springer-Verlag, 1981.
\bibitem{MW} I.~Madsen, M.~S.~Weiss, The stable moduli space of Riemann surfaces: Mumford's conjecture, Ann. of Math., 165 (2007), 843--941.
\bibitem{EM0} N.~Mishachev, Y.~Eliashberg, Surgery of singularities of foliations, Funct. Anal. and Pril., 11 (1977), 43--53; English translation in Funct. Anal. Appl., 11 (1977), 197--206. 
\bibitem{Ph} A.~Phillips, Submersions of open manifolds, Topol., 6 (1967), 171--206.
\bibitem{Pl} A.~du Plessis, Contact-invariant regularity
conditions, Lect. Notes Math., 535, Springer-Verlag, 1976,
205--236.
\bibitem{Pl1} A.~du Plessis, Maps without certain singularities, Comment. Math. Helv., 50(1975), 363--382.
\bibitem{Ri} R.~Rim\'{a}nyi, Thom polynomials, symmetries and incidences of singularities, Inv. Math., 143 (2001), 499--521.
\bibitem{RS} R.~Rim\'{a}nyi, A.~Sz\H{u}cs, Pontrjagin-Thom type
construction for maps with singularities, Topology, 37 (1998),
1177--1191.
\bibitem{RoS}  C.~Rourke, B.~Sanderson, The compression theorem. I, Geom. Topol., 5 (2001), 399--429. 
\bibitem{Ru} Y.~Rudyak, On Thom spectra, orientability, and cobordism, Springer monographs in math., Springer Verlag, 1998.
\bibitem{Sad1} R.~Sadykov, Bordism groups of special generic
mappings, Proc. Amer. Math. Soc., 133 (2005), 931--936.
\bibitem{Sad5} R.~Sadykov, Bordism groups of solutions to differential relations, arXiv:math/0608460. 
\bibitem{Sad3} R.~Sadykov, Singular cobordism categories, preprint, arXiv/math:0804.1267, 2008. 
\bibitem{Sad4} R.~Sadykov, Cobordism groups of Morin maps, preprint, 2008. 
\bibitem{Sae1} O.~Saeki, Cobordism groups of special generic functions and groups of homotopy spheres, Japan. J. Math., 28 (2002), 287--297.
\bibitem{SY} O.~Saeki, T.~Yamamoto, Singular fibers and characteristic classes, Topology Appl. 155 (2007), 112--120.
\bibitem{Sp} D.~Spring, Convex integration theory: solutions to
the h-principle in geometry and topology, Birkh\"auser Verlag,
Basel, Boston, Berlin, 1998.
\bibitem{St} E.~Stong, Notes on cobordism theory, Princeton Univ. Press and Univ. Tokyo Press, 1968.
\bibitem{Sw} R.~M.~Switzer, Algebraic topology--homology and homotopy, Springer, 1975. 
\bibitem{Sz} A.~Sz\H ucs, Elimination of singularities by
cobordism, Real and Complex Singularities (S\~ao Carlos, 2002), Ed. T.Gaffney and M.Ruas, Contemp. Math., 354, AMS, Providence, RI (2004), 301--324.
\bibitem{Sz1} A.~Sz\H ucs, Cobordism of singular maps, Geom. and Topol., 12 (2008), 2379--2452. 
\bibitem{Th1} W.~P.~Thurston, The theory of foliations of codimension greater than one, Comm. Math. Helvet., 49 (1974), 214--231. 
\bibitem{Th2} W.~P.~Thurston, Existence of codimension-one foliations, Ann. of Math., 104 (1976), 249--268.
\bibitem{Wa} C.~T.~C.~Wall, A second note on symmetry of singularities, Bull. London Math. Soc., 12 (1980), 347--354.
\bibitem{We} R.~Wells, Cobordism groups of immersions, Topology, 5
(1966), 281--294.
\end{thebibliography}
\end{document}